\documentclass[]{interact}
\usepackage{epstopdf}

\usepackage{graphicx}
\usepackage{algorithm}
\usepackage{algpseudocode}
\usepackage{amssymb}
\usepackage{amsmath}
\usepackage{subcaption}
\usepackage{xcolor}
\usepackage{breqn}
\usepackage{booktabs} 
\usepackage{multirow}

\usepackage[hidelinks]{hyperref}
\usepackage{amsthm}
\newtheorem{assumption}{Assumption}
\newtheorem{proposition}{Proposition}
\newtheorem{remark}{Remark}

\newcommand{\vf}{\mathbf{f}}
\newcommand{\vx}{\mathbf{x}}
\newcommand{\vu}{\mathbf{u}}
\newcommand{\vh}{\mathbf{h}}
\newcommand{\vU}{\mathbf{U}}

\newcommand{\vZ}{\mathbf{Z}}
\newcommand{\vz}{\mathbf{z}}
\newcommand{\vS}{\mathbf{S}}
\newcommand{\bv}{\mathbf{v}}
\newcommand{\vl}{\mathbf{l}}
\newcommand{\vQ}{\mathbf{Q}}
\newcommand{\vR}{\mathbf{R}}
\newcommand{\vk}{\mathbf{k}}
\newcommand{\vK}{\mathbf{K}}
\newcommand{\blambda}{\boldsymbol{\lambda}}
\newcommand{\bLambda}{\boldsymbol{\Lambda}}

\newcommand{\bo}{\boldsymbol{0}}
\newcommand{\cJ}{\mathcal{J}}
\newcommand{\tran}{^{\top}}
\makeatletter
\long\def\@maketablecaption#1#2{\@tablecaptionsize
  \setbox\@tempboxa\hbox{#1. #2}
  \ifdim \wd\@tempboxa >\hsize              
    \unhbox\@tempboxa\par                   
  \else                                     
    \global \@minipagefalse
    \hbox to\hsize{\hfil\box\@tempboxa\hfil}
  \fi}
\makeatother

\begin{document}


\title{A Sequential Quadratic Programming Perspective on Optimal Control}

\author{
\name{
Abhijeet\thanks{Email: abhinir@tamu.edu}
and
Suman Chakravorty\thanks{Email: schakrav@tamu.edu}
}
\affil{Texas A\&M University, College Station, TX, USA}
}

\maketitle

\begin{abstract}                        This paper investigates the relationship between the iterative Linear Quadratic Regulator (iLQR), Differential Dynamic Programming (DDP), and Newton's method for optimal control problems. Using a unified Sequential Quadratic Programming (SQP) framework, we introduce Newton-LQR, a Riccati-based algorithm that solves for the exact Newton step the optimal control problem. Within this framework, we establish that iLQR is a principled modified SQP method rather than merely an approximation of DDP. We show that iLQR generates a descent direction, whereas the inclusion of dynamics Hessians in DDP and Newton-LQR may lead to indefinite quadratic subproblems and necessitate regularization. We further distinguish DDP and the proposed Newton-LQR algorithm through their treatment of costate information and second-order dynamics terms. Numerical experiments across six benchmark problems demonstrate that iLQR often exhibits superior convergence during the initial optimization phase, while Newton-LQR and DDP can exploit second-order information to accelerate convergence near a solution. 
\end{abstract}

\section{Introduction}

The study of iterative trajectory optimization has a long history in nonlinear
control. Among the most influential methods are Differential Dynamic Programming
(DDP), introduced in the 1960s \cite{mayne1965ddp, JacobsonMayne1970DDP}, and the
iterative Linear Quadratic Regulator (iLQR), which later emerged as a simplified
yet powerful alternative \cite{todorov2005ilqg, tassa2012synthesis}. DDP is acclaimed as
Newton-type algorithm for trajectory optimization that exploits second-order
expansions of both the cost and the dynamics, achieving locally quadratic
convergence near an optimum \cite{mayne1965ddp, JacobsonMayne1970DDP,
jacobson1970ddp}. iLQR, by contrast, employs consistent first-order
approximations of the dynamics, making it more numerically reliable and
computationally efficient in practice \cite{todorov2005ilqg, tassa2012synthesis,D2C}.
Over the decades, these two approaches have become foundational tools for
nonlinear control, and trajectory optimization
\cite{BrysonHo69,Wang2025SearchFeedbackRL,inf_horizon1,inf_horizon2}.

A large body of literature has analyzed the theoretical properties and limitations
of both methods. Several works have emphasized the quadratic convergence of DDP
under ideal assumptions \cite{mayne1965ddp, JacobsonMayne1970DDP, Pantoja_SN, shoemaker1983ddp},
but have also noted the difficulties that arise in practice when the algorithm's
backward pass combines second-order expansions with first-order terms, which can
compromise stability and necessitate heavy regularization \cite{Pantoja_SN,LiaoShoemaker1992,feng2020optimization,NgangaWensing2021RAL,pardo2017time}.
An explicit connection between DDP and Newton's method was established in
\cite{Pantoja_SN,Murray1984DDP_Newton}, and it is argued in \cite{LiaoShoemaker1992} that, for
discrete-time optimal control, DDP can outperform Newton's method---retaining
strong local convergence while offering greater robustness and efficiency under
practical problem structures. For rigid-body systems, accelerated second-order DDP
has been explored \cite{NgangaWensing2021RAL}. Extensions such as the iterative
Linear Quadratic Gaussian (iLQG) \cite{todorov2005ilqg} and control-limited iLQR
\cite{tassa2014control} have reinforced the practical importance of iLQR in
stochastic and constrained settings.



Recent developments in iLQR and DDP have significantly expanded their applications in nonlinear trajectory optimization, particularly in robotics and aerospace engineering. In robotics, these methods have enabled efficient motion planning and control for high-dimensional dynamical systems. Notable algorithmic developments include ALTRO, which combines augmented Lagrangian methods with iLQR-based trajectory optimization \cite{Howell2019ALTRO}, and Crocoddyl, which introduces a feasibility-driven DDP framework for multi-contact optimal control \cite{Mastalli2020Crocoddyl}. Further advances by Wensing and collaborators have extended DDP to hybrid dynamical systems for whole-body motion planning of legged robots \cite{LiWensing2020HybridDDP}, developed computationally efficient approaches for incorporating second-order dynamics information into DDP \cite{NgangaWensing2021Accelerating}, and demonstrated full second-order DDP for whole-body model predictive control on legged robot hardware \cite{NgangaLiWensing2023MPC}. In aerospace engineering, DDP has been applied to trajectory optimization in astrodynamics and spacecraft guidance \cite{Russell1,Russell2,Aziz2019HDDP_CR3BP}, as well as planetary landing problems \cite{Noyes2021RobustEntryGuidance}. Other extensions, such as parameterized DDP, have demonstrated its applicability to complex robotic systems \cite{Oshin2022PDDP_RSS}. While effective in these contexts, these studies also highlight the algorithm's sensitivity to problem structure and its dependence on strong regularization, particularly when the local quadratic models become inconsistent \cite{Pantoja_SN,LiaoShoemaker1992}.

Although iLQR is often described as an approximation of DDP obtained by
neglecting second-order dynamic terms \cite{todorov2005ilqg}, this characterization
obscures its principled nature. In this paper, we revisit the solution of optimal
control problems through the lens of  Sequential Quadratic Programming
(SQP) \cite{boggs1995sqp, Nocedal2006NumericalOptimization}.The critical insight is that 
the Newton step in a constrained optimization problem is equivalent to solving a
Quadratic Program (QP), which reduces to an LQR problem in the context of optimal
control. Using this framework, we show that iLQR is not merely an approximation
of DDP, but rather a principled method that guarantees cost-descent directions
under mild assumptions---a property that neither DDP nor Newton's method can ensure
in general, particularly when far from an optimum \cite{Pantoja_SN,
LiaoShoemaker1992}. We also derive the exact Newton LQR formulation, and
show that DDP can be interpreted as an approximation of this exact Newton step
only near an optimum.To evaluate these theoretical findings, we conduct numerical comparisons on pendulum and cart-pole swing-up tasks, along with four additional nonlinear benchmark problems. Our analysis examines the positive definiteness of the control Hessian, the frequency and magnitude of regularization required to obtain descent directions, and the resulting convergence behavior. The numerical results corroborate the SQP-based analysis, demonstrating that iLQR generally maintains positive-definite quadratic subproblems, requires substantially less regularization, and achieves faster convergence than DDP across the problems considered. Furthermore, by avoiding the explicit computation of dynamics Hessians, iLQR reduces the computational overhead associated with constructing the local quadratic models \cite{Wang2025SearchFeedbackRL,D2C}. We also see that close to an optimum, Newton LQR and DDP exhibit faster convergence than iLQR. These observations highlight the practical advantages of iLQR, particularly during the initial stages of optimization when the iterates are far from a stationary solution, while also motivating the switch to second-order dynamics information in Newton-LQR and DDP to accelerate convergence near optimality.

The remainder of this paper is organized as follows. Section~\ref{sec:background} introduces the optimal control problem and establishes the necessary notation. Section~\ref{sec:sqp} presents the Sequential Quadratic Programming (SQP) framework. Section~\ref{sec:ipfocs} develops the Newton-LQR formulation and examines its relationship with iLQR, DDP, and Stagewise-Newton. Section~\ref{sec:algorithms} presents numerical comparisons across several benchmark problems and investigates the proposed hybrid approach. Finally, Section~\ref{sec:conclusion} summarizes the main findings and conclusions.
Section~\ref{sec:conclusion} concludes with insights drawn from these comparisons.

\section{Background}\label{sec:background}
We address the general class of discrete-time, finite-horizon optimal control problems. Let the state of a system at time step $t$ be $\vx_t \in \mathbb{R}^n$, and the control input be $\vu_t \in \mathbb{R}^m$. The system evolves over a finite time horizon of $T$ steps according to the nonlinear state-transition function $f$:
\begin{equation}
    \vx_{t+1} = \vf(\vx_t, \vu_t), \qquad t=0, 1, \cdots, T-1,
    \label{eq:dynamics}
\end{equation}
starting from a given initial state $\vx_0$.

The objective is to find a sequence of control inputs $\mathbf{U} = \{\vu_0, \dots, \vu_{T-1}\}$ that minimizes a total cost functional, $\cJ$. The optimal control problem can be formally stated as:
\begin{equation}
    \min_{\mathbf{U}} \cJ(\mathbf{U}) = \sum_{t=0}^{T-1} c(\vx_t, \vu_t) + C_T(\vx_T),
    \label{eq:ocp}
\end{equation}
subject to the dynamics in \eqref{eq:dynamics}. To ensure the problem is well-posed for the optimization methods used in this work, we make the following assumption.

\begin{assumption}
The system dynamics function $\vf$, stage cost $c$, and terminal cost $C_T$ are twice continuously differentiable ($\mathcal{C}^2$) functions with respect to their arguments.
\label{assum:smoothness}
\end{assumption}
We introduce the Newton-LQR approach to solve the optimal control problem described in eq.~\eqref{eq:ocp}. Along with this, we analyze the prominent trajectory optimization methods, including Differential Dynamic Programming (DDP), the iterative Linear Quadratic Regulator (iLQR), and Stagewise Newton methods, all rely on a common underlying principle. They iteratively refine a nominal trajectory by solving a local, second-order approximation of the original optimal control problem. Constructing this approximation requires the derivatives of the dynamics and cost functions evaluated along the nominal trajectory.

Let $(\bar{\vx}_t, \bar{\vu}_t)$ for $t=0, \dots, T-1$ denote such a nominal trajectory, which is a dynamically feasible sequence satisfying \eqref{eq:dynamics}. We define the state and control perturbations, or variations, from this nominal trajectory as:
\begin{equation}
    \delta \vx_t := \vx_t - \bar{\vx}_t \quad \text{and} \quad \delta \vu_t := \vu_t - \bar{\vu}_t.
\end{equation}

For notational brevity, we adopt a compact subscript notation for the partial derivatives. The first-order derivatives (Jacobians) and second-order derivatives (Hessians) of the dynamics and cost functions are written as:
\[
\vf_{\vx_t}, \quad \vf_{\vu_t}, \quad \vf_{\vx_t\vx_t}, \quad \vf_{\vu_t\vu_t}, \quad \vf_{\vu_t\vx_t}, \quad
c_{\vx_t}, \quad c_{\vu_t}, \quad c_{\vx_t\vx_t}, \quad c_{\vu_t\vu_t}, \quad c_{\vx_t\vu_t}.
\]
Derivatives of the terminal cost are denoted similarly, e.g., $C_{T,\vx}$ and $C_{T,\vx\vx}$. Unless stated otherwise, all such derivatives are evaluated at the corresponding time step along the nominal trajectory, i.e., $\vf_\vx \equiv \frac{\partial \vf}{\partial \vx}(\bar{\vx}_t, \bar{\vu}_t)$. We use $\otimes$ to denote a tensor product. Finally, we denote the transpose of a vector or matrix $\mathbf{A}$ as $\mathbf{A}\tran$.

A unified perspective for understanding and comparing the trajectory optimization algorithms in this paper is offered by the framework of Sequential Quadratic Programming (SQP). We therefore begin by developing the fundamental concepts of the SQP approach as applied to the general optimal control problem described in eq.~\eqref{eq:ocp}. 

\section{Sequential Quadratic Programming}\label{sec:sqp}
In this section, we outline the basic SQP program for future references. Consider the generic optimization problem:
\begin{align*}
    &\min_{\vz} g(\mathbf{z}),\\
    \text{subject to:}~~&\mathbf{h}(\mathbf{z}) = \boldsymbol{0}.
\end{align*}
This constrained optimization can be converted to an unconstrained optimization using an augmented Lagrangian as:
\begin{align*}
    g(\mathbf{z}) + \bLambda\tran \mathbf{h}(\mathbf{z}) = 0,
\end{align*}
where $\bLambda$ denotes the set of Lagrange multipliers. This gives the necessary conditions as:
\[
\frac{\partial \mathcal{L}}{\partial \mathbf{z}} = \boldsymbol{0} 
\quad \text{and} \quad 
\frac{\partial \mathcal{L}}{\partial \boldsymbol{\bLambda}} = \boldsymbol{0},
\]
yielding:
\[
g_{\mathbf{z}} +  \vh_{\mathbf{z}}\tran \bLambda^*= \boldsymbol{0},\;
\vh(\mathbf{z}^*) = \boldsymbol{0},
\]
at a stationary point \((\vz^*, \bLambda^*)\).

To solve the problem, the second-order Newton method suggests to start with a nominal trajectory: \((\bar{\vz}, \bar{\bLambda})\) and then expand the Lagrangian 
\(\mathcal{L}\) to second order around the nominal.
\begin{align*}
\mathcal{L}(\vz, \bLambda) = g + \bLambda \tran \mathbf{h} 
= \big( \bar{g} + \bar{\bLambda}\tran\bar{\mathbf{h}} \big) 
+ (g_\vz \delta \vz + \bar{\bLambda}\tran \mathbf{h}_\vz \delta \vz + \delta \bLambda\tran\bar{\mathbf{h}}) \\
+ \frac{1}{2} \delta \vz\tran \big( g_{\vz  \vz} + \bar{\bLambda} \otimes \mathbf{h}_{\vz \vz} \big) \delta \vz 
+  \delta \vz\tran \, \mathbf{h}_\vz \delta \bLambda 
 \\
\equiv \bar{\mathcal{L}} + \Delta \mathcal{L},
\end{align*}
where $\bar{g} = g(\bar{\vz})$, $\bar{\mathbf{h}} = \mathbf{h}(\bar{\vz})$, $\bar{\mathcal{L}} = \bar{g}+\bar{\bLambda}\tran\bar{\mathbf{h}}$.


In a standard SQP approach, one follows iterative Newton's method and optimizes the \((\Delta \mathcal{L})\) part in terms of the perturbations \((\delta \vz, \delta \bLambda)\) and updates the solution iteratively as:
\[
\vz = \bar{\vz} + \delta \vz, 
\quad 
\bLambda = \bar{\bLambda} + \delta \bLambda.
\]

The necessary conditions for obtaining the optimum of \(\Delta \mathcal{L}\) in terms of \(\delta x, \delta \lambda\) are obtained by setting
\[
\frac{\partial (\Delta \mathcal{L})}{\partial (\delta \vz)} = 0 
\quad \text{and} \quad 
\frac{\partial (\Delta \mathcal{L})}{\partial (\delta \bLambda)} = 0,
\]
yielding:
\[
\begin{bmatrix}
g_{\vz \vz} + \bar{\bLambda}\otimes \mathbf{h}_{\vz \vz} && \mathbf{h}_{\vz}\tran \\
\mathbf{h}_{\vz} && \boldsymbol{0}
\end{bmatrix}
\begin{bmatrix}
\delta \vz \\[6pt]
\delta \bLambda
\end{bmatrix}
=
\begin{bmatrix}
- g_{\vz} - \mathbf{h}_{\vz}\tran\bar{\bLambda} \\[6pt]
-\bar{\mathbf{h}}
\end{bmatrix}.
\]
However, taking the \( \mathbf{h}_\vz\tran \bar{\bLambda}\) term to the LHS gives:
\[
\begin{bmatrix}
g_{\vz \vz} + \bar{\bLambda}\otimes \mathbf{h}_{\vz \vz} && \mathbf{h}_{\vz}\tran \\
\mathbf{h}_{\vz} && \boldsymbol{0}
\end{bmatrix}
\begin{bmatrix}
\delta \vz \\[6pt]
\bar{\bLambda} + \delta \bLambda
\end{bmatrix}
=
\begin{bmatrix}
- g_{\vz} \\[6pt]
-\bar{\mathbf{h}}
\end{bmatrix},
\]
where \(\bLambda = \bar{\bLambda} + \delta \bLambda\) is simply the new \(\bLambda\) (the updated $\bLambda$ at new iteration).

However, the necessary conditions above are simply the solution to the QP problem:
\begin{align}
\min_{\delta \vz} \; g_\vz \delta \vz + \tfrac{1}{2} \delta \vz\tran \big(g_{\vz\vz} + \bar{\bLambda} \otimes \mathbf{h}_{\vz \vz}\big) \delta \vz \nonumber \\
\text{s.t.} \quad \mathbf{h}_{\vz} \delta \vz + \bar{\mathbf{h}} = 0. \label{N-QP}
\end{align}

Hence, the Newton step can be accomplished by solving the above QP for the solution \((\delta \vz^q, \bLambda^q)\) and updating the solution as:
\[
\vz = \bar{\vz} + \delta \vz^q, 
\quad 
\bLambda = \bLambda^q.
\]

Thus, the original constrained optimization can be solved by solving a sequence of QP problems above, and this will be termed as the Newton SQP approach. \\
However, in general, it is not advisable to use the ``true" Hessian of the Lagrangian \(\big(g_{\vz \vz} + \bar{\bLambda}\otimes \mathbf{h}_{\vz \vz}\big)\) in the QP unless one is close to an optimum.

To see this, note that in the constrained case, the descent direction 
$\delta \vZ = \begin{pmatrix} \delta \vz \\ \delta \bLambda \end{pmatrix}$ 
is given by:
\[
\delta \vZ = -(\nabla^2 \mathcal{L})^{-1} \nabla \mathcal{L},
\]
and
\[
\nabla^2 \mathcal{L} = 
\begin{bmatrix}
g_{\vz \vz} + \bar{\bLambda} \otimes \mathbf{h}_{\vz \vz} && \mathbf{h}_{\vz}\tran \\
\mathbf{h}_{\vz} && \boldsymbol{0}
\end{bmatrix}
\quad \nsucc 0
\]
if the current solution is far from an optimum.  

If the Hessian $\nabla^2 \mathcal{L} \succ 0$, the Newton direction can be assured to be a descent direction, since then:
\[
\delta \vZ\tran \nabla \mathcal{L} = -\nabla \mathcal{L}\tran (\nabla^2 \mathcal{L})^{-1} \nabla \mathcal{L} < 0.
\]
However, if $\nabla^2 \mathcal{L} \nsucc 0$, there is no such guarantee.  

Thus, in order to be able to use the SQP approach, it becomes critical that one be able to ensure that the QP solution is a descent direction for the cost function.  

Consider now a modified problem (QP):
\begin{align}
\min_{\delta \vz} \; g_\vz \delta \vz + \tfrac{1}{2} \, \delta \vz\tran g_{\vz \vz} \delta \vz \nonumber\\
\mathbf{h}_{\vz} \delta \vz + \bar{\mathbf{h}} = 0. \label{M-QP}
\end{align}

We make the assumption that $\mathbf{h}(\bar{\vz}) = \bar{\mathbf{h}} = \boldsymbol{0}$, i.e., our current estimate is feasible.  

The necessary conditions for this QP are:
\begin{align}
g_\vz + g_{\vz \vz} \delta \vz + \mathbf{h}_{\vz}\tran \bLambda &= \boldsymbol{0},  \\
\bar{\mathbf{h}} + \mathbf{h}_\vz \delta \vz &= \boldsymbol{0}. 
\end{align}



Since $g(\vz)$ is a cost function by design, it is a reasonable assumption that $g_{\vz \vz} \succ 0$.  

Hence, the solution to the modified QP \eqref{M-QP} is always a descent direction for the cost function.  
Furthermore,
\[
\mathbf{h}(\bar{\vz} + \delta \vz) \approx \mathbf{h}(\bar{\vz}) + \mathbf{h}_\vz \delta \vz \approx \boldsymbol{0}
\]
if $\delta \vz$ is small, and thus the update $\vz = \bar{\vz} + \delta \vz$ also satisfies the constraint. If $\delta \vz$ is not small enough, we can do a line search using a parameter $\alpha \in (0,1]$ such that
\[
\mathbf{h}(\vz) = \mathbf{h}(\bar{\vz} + \alpha \delta \vz) \approx \mathbf{h}(\bar{\vz}) + \alpha \mathbf{h}_\vz \delta \vz \approx 0,
\]
and thus $\vz = \bar{\vz} + \alpha \delta \vz$ is feasible.  

Hence, if we have a feasible solution, the modified QP always gives a descent direction for the cost while ensuring that the solution remains feasible if allied with a suitable line search approach.  
The development above can be summarized as the following result:  

\medskip

\begin{proposition} \label{prop:mod_QP}
Consider the modified QP \eqref{M-QP}. Given $g_{\vz \vz} \succ 0$, the solution to the QP, $\delta \vz$, is always a descent direction for the cost function $g(\vz)$. Furthermore, the new solution $\vz = \bar{\vz} + \alpha \delta \vz$ can be made feasible, i.e., $\mathbf{h}(\vz) = \boldsymbol{0}$, by a suitable choice of the line-search parameter, $\alpha \in (0,1]$.
\end{proposition}

\section{The SQP perspective on iterative perturbation feedback methods}\label{sec:ipfocs}
In Section~\ref{sec:sqp}, we introduced the exact Newton
QP~\eqref{N-QP} and its modified counterpart~\eqref{M-QP}.
We now specialize these formulations to the optimal
control problem~\eqref{eq:ocp} and demonstrate how they
give rise to the Newton-LQR and iLQR algorithms,
respectively. We subsequently examine DDP and
Stagewise-Newton as related approximations to the
Newton-LQR formulation. This unified SQP perspective
allows us to identify the differences in the quadratic
subproblems solved by these algorithms and their
implications for convergence. These methods, which we collectively refer to as ``iterative perturbation feedback methods,'' exploit the sparsity structure of the optimal control problem and are intended to solve the problem defined by eq.~\eqref{eq:ocp} and eq.~\eqref{eq:dynamics}.

All of these algorithms are initialized from a nominal trajectory, $(\bar{\vx}_t,\bar{\vu}_t)$, which is iteratively refined until convergence. To clarify the relationship among them, we first revisit the QP subproblem obtained from the exact Newton linear-quadratic approximation. This serves as the natural reference point for the discussion that follows. We then consider the modified approximation underlying iLQR, and finally examine the formulations used by DDP and Stagewise-Newton. As will become clear, although DDP and Stagewise-Newton are often interpreted as Newton-type methods, the QPs they solve differ from the exact Newton subproblem described in the next subsection and have to be understood as approximations.
\subsection{The Newton Step in Optimal Control}
We begin with the following LQR problem, obtained by expanding the optimal control problem about the nominal trajectory:
\begin{align}
    &\min_{\delta \vu_{t}} \sum_{t=0}^{T-1} \Bigg( c_{\vx_t}\delta \vx_{t} + c_{\vu_t}\delta \vu_{t} + \frac{1}{2} \delta \vx_t\tran c_{\vx_t\vx_t} \delta \vx_t + \delta \vu_t\tran c_{\vu_t\vx_t} \delta \vx_t + \frac{1}{2} \delta \vu_t\tran c_{\vu_t\vu_t} \delta \vu_t \nonumber \\
    &\quad \quad \quad + \frac{\bar{\blambda}_{t+1}}{2}\begin{bmatrix}
        \delta \vx_{t}\\
        \delta \vu_{t}
    \end{bmatrix}\tran \begin{bmatrix}
        \vf_{\vx_t\vx_t} && \vf_{\vx_t\vu_t}\\
        \vf_{\vu_t\vx_t} && \vf_{\vu_t\vu_t}
    \end{bmatrix}\begin{bmatrix}
        \delta \vx_{t}\\
        \delta \vu_{t}
    \end{bmatrix} \Bigg) + C_{T,\vx}\delta \vx_{T} + \frac{1}{2} \delta \vx_{T}\tran C_{T,\vx\vx} \delta \vx_{T},\nonumber \\
    &\text{subject to:}~\delta \vx_{t+1} = \vf_{\vx_t}\delta \vx_{t} + \vf_{\vu_t}\delta \vu_{t},\nonumber\\
    &\text{given:}~~\delta \vx_{0} = \bo.\label{eq:DDP-LQR}
\end{align}
The term $\bar{\boldsymbol{\lambda}}_{t}$ represents the Lagrange multipliers from the previous iteration. The above problem is obtained by applying the Newton QP described in the previous section in eq.~\eqref{N-QP} to the optimal control problem. Essentially, the following result shows the state and control update, along with the new set of Lagrange multipliers which are then used recursively in a Newton scheme. Specifically, the quadratic objective corresponds to
the second-order expansion of the Lagrangian in
\eqref{N-QP}, while the equality constraints correspond
to the linearized system dynamics. Consequently,
the solution of this LQR problem yields the exact
Newton/SQP direction for a dynamically feasible
nominal trajectory.
\begin{proposition}
The solution of the Newton--LQR problem in~\eqref{eq:DDP-LQR} is characterized by the forward pass:
\begin{align*}
    \delta \vx_{t+1} &= \vf_{\vx_t}\,\delta \vx_t + \vf_{\vu_t}\,\delta \vu_t,\\
    \blambda_t &= \vS_t \,\delta \vx_t + \bv_t,\\
    \delta \vu_t &= \vk_t + \vK_t \,\delta \vx_t,
\end{align*}
where the feedforward and feedback gains are given by
\[
\vk_t = -(Q_{\vu\vu}^t)^{-1}Q_{\vu}^t,
\qquad
\vK_t = -(Q_{\vu\vu}^t)^{-1}Q_{\vu\vx}^t,
\]
and they are computed in the backward pass along with $\vS_{t}$ and $\bv_{t}$ as
\begin{align}
    \vS_t &= Q_{\vx\vx}^t - Q_{\vx\vu}^t (Q_{\vu\vu}^t)^{-1} Q_{\vu\vx}^t
    \qquad \text{with} \qquad
    \vS_{T} = C_{T,\vx\vx}, \label{eq:back_pass_S}\\
    \bv_t &= Q_{\vx}^t - Q_{\vx\vu}^t (Q_{\vu\vu}^t)^{-1} Q_{\vu}^t
    \qquad \text{with} \qquad
    \bv_{T} = C_{T,\vx}. \label{eq:back_pass_v}
\end{align}
The corresponding coefficients are
\begin{align*}
    Q_{\vx}^t
    &= c_{\vx_t} + \vf_{\vx_t}\tran \bv_{t+1},\\
    Q_{\vu}^t
    &= c_{\vu_t} + \vf_{\vu_t}\tran \bv_{t+1},\\
    Q_{\vx\vx}^t
    &= c_{\vx_t\vx_t}
    + \vf_{\vx_t}\tran \vS_{t+1}\vf_{\vx_t}
    + \bar{\blambda}_{t+1} \otimes \vf_{\vx_t\vx_t},\\
    Q_{\vx\vu}^t
    &= c_{\vx_t\vu_t}
    + \vf_{\vx_t}\tran \vS_{t+1}\vf_{\vu_t}
    + \bar{\blambda}_{t+1} \otimes \vf_{\vx_t\vu_t},\\
    Q_{\vu\vu}^t
    &= c_{\vu_t\vu_t}
    + \vf_{\vu_t}\tran \vS_{t+1}\vf_{\vu_t}
    + \bar{\blambda}_{t+1} \otimes \vf_{\vu_t\vu_t}.
\end{align*}
\end{proposition}

\begin{proof}
Consider the Hamiltonian
\begin{align*}
    H =  &c_{\vx_t}\delta \vx_{t} + c_{\vu_t}\delta \vu_{t} + \frac{1}{2} \delta \vx_t\tran c_{\vx_t\vx_t} \delta \vx_t + \delta \vu_t\tran c_{\vu_t\vx_t} \delta \vx_t + \frac{1}{2} \delta \vu_t\tran c_{\vu_t\vu_t} \delta \vu_t \\
    &+ \frac{\bar{\blambda}_{t+1}}{2}\begin{bmatrix}
        \delta \vx_{t}\\
        \delta \vu_{t}
    \end{bmatrix}\tran \begin{bmatrix}
        \vf_{\vx_t\vx_t} && \vf_{\vx_t\vu_t}\\
        \vf_{\vu_t\vx_t} && \vf_{\vu_t\vu_t}
    \end{bmatrix}\begin{bmatrix}
        \delta \vx_{t}\\
        \delta \vu_{t}
    \end{bmatrix}
    + \blambda_{t+1}\tran \big(\vf_{\vx_t}\delta \vx_{t} + \vf_{\vu_t}\delta \vu_{t} - \delta \vx_{t+1} \big).
\end{align*}
The necessary conditions are
\begin{align*}
    \frac{\partial H}{\partial \vu_{t}} = 0,
    \qquad
    \frac{\partial H}{\partial \blambda_{t}} = 0,
    \qquad \text{and} \qquad
    \frac{\partial H}{\partial \vx_{t}} = \blambda_{t}.
\end{align*}
The transversality condition is
\begin{align*}
    \blambda_{T} = C_{T,\vx} + C_{T,\vx\vx}\delta \vx_{T}.
\end{align*}

Owing to the LQR structure of the problem, from induction, the costate have the affine form:
\begin{align}
    \blambda_{t} = \vS_{t}\delta \vx_{t} + \bv_{t}.
    \label{eq:lambda_form}
\end{align}
Using the condition $\frac{\partial H}{\partial \blambda_t}=0$, we obtain
\begin{align}
    \vf_{\vx_t} \delta \vx_{t} + \vf_{\vu_t} \delta \vu_t - \delta \vx_{t+1} = 0,
    \label{eq:per_dyn_lin_v1}
\end{align}
or equivalently,
\[
\delta \vx_{t+1} = \vf_{\vx_t} \delta \vx_{t} + \vf_{\vu_t} \delta \vu_t.
\]
Next, from the necessary condition, $\frac{\partial H}{\partial \vu_t}=0$, we have
\begin{align*}
    c_{\vu_t} +& c_{\vu_t\vu_t}\delta \vu_t + c_{\vu_t\vx_t}\delta \vx_t\\
    &+ (\bar{\blambda}_{t+1} \otimes \vf_{\vu_t\vx_t})\delta \vx_t
    + (\bar{\blambda}_{t+1} \otimes \vf_{\vu_t\vu_t})\delta \vu_t
    + \vf_{\vu_t}\tran\blambda_{t+1} = 0.
\end{align*}
Substituting the affine representation~\eqref{eq:lambda_form} for $\blambda_{t+1}$ yields
\begin{align*}
    0 =\;& c_{\vu_t} + c_{\vu_t\vu_t}\delta \vu_t + c_{\vu_t\vx_t}\delta \vx_t
    + (\bar{\blambda}_{t+1} \otimes \vf_{\vu_t\vx_t})\delta \vx_t
    + (\bar{\blambda}_{t+1} \otimes \vf_{\vu_t\vu_t})\delta \vu_t \\
    &+ \vf_{\vu_t}\tran(\vS_{t+1}\delta \vx_{t+1}+\bv_{t+1}).
\end{align*}
Using~\eqref{eq:per_dyn_lin_v1}, this becomes
\begin{align*}
    0 =\;& c_{\vu_t} + c_{\vu_t\vu_t}\delta \vu_t + c_{\vu_t\vx_t}\delta \vx_t
    + (\bar{\blambda}_{t+1} \otimes \vf_{\vu_t\vx_t})\delta \vx_t
    + (\bar{\blambda}_{t+1} \otimes \vf_{\vu_t\vu_t})\delta \vu_t \\
    &+ \vf_{\vu_t}\tran\!\big(\vS_{t+1}(\vf_{\vx_t}\delta \vx_t + \vf_{\vu_t}\delta \vu_t)+\bv_{t+1}\big).
\end{align*}
Collecting terms in $\delta \vu_t$ and $\delta \vx_t$, we obtain
\begin{align*}
    \delta \vu_{t} = &-\left( c_{\vu_t\vu_t} + \bar{\blambda}_{t+1} \otimes \vf_{\vu_t\vu_t} + \vf_{\vu_t}\tran \vS_{t+1}\vf_{\vu_{t}}\right)^{-1}\\
    &\cdot \left( \left( c_{\vu_t} + \vf_{\vu_t}\tran \bv_{t+1}  \right)
    + \left( c_{\vu_t\vx_t} + \bar{\blambda}_{t+1} \otimes \vf_{\vu_t\vx_t} + \vf_{\vu_t}\tran \vS_{t+1}\vf_{\vx_{t}}\right) \delta \vx_{t} \right).
\end{align*}
The above equations can be written in terms of $Q_{\vu\vu}^{t},Q_{\vu}^{t},~\text{and}~Q_{\vu\vx}^{t}$ as:
\begin{align}
    \delta \vu_{t} = -(Q_{\vu\vu}^{t})^{-1}Q_{\vu}^{t} -(Q_{\vu\vu}^{t})^{-1}Q_{\vu\vx}^{t}\delta \vx_{t}
    = \vk_{t} + \vK_{t}\delta \vx_{t}.
    \label{eq:control_update_v1}
\end{align}

We now turn to the costate equation. From the necessary condition $\frac{\partial H}{\partial \vx_t}=\blambda_t$, we obtain
\begin{align*}
    \blambda_{t} =\;& c_{\vx_t} + c_{\vx_t\vx_t}\delta \vx_{t} + c_{\vx_t\vu_t} \delta \vu_t
    + (\bar{\blambda}_{t+1}\otimes \vf_{\vx_t\vx_t})\delta \vx_{t} \\
    &+ (\bar{\blambda}_{t+1}\otimes \vf_{\vx_t\vu_t})\delta \vu_{t}
    + \vf_{\vx_t}\tran\blambda_{t+1}.
\end{align*}
Substituting the affine form~\eqref{eq:lambda_form} for $\blambda_t$ and $\blambda_{t+1}$, together with the state recursion~\eqref{eq:per_dyn_lin_v1}, gives
\begin{align*}
    \vS_{t}\delta \vx_{t} + \bv_{t} =\;& c_{\vx_t} + c_{\vx_t\vx_t}\delta \vx_{t}
    + \left( c_{\vx_t\vu_t}+\bar{\blambda}_{t+1}\otimes \vf_{\vx_t\vu_t} \right)  \delta \vu_{t} \\
    &+ (\bar{\blambda}_{t+1}\otimes \vf_{\vx_t\vx_t})\delta \vx_{t}
    + \vf_{\vx_t}\tran\!\left(\vS_{t+1}\left( \vf_{\vx_t}\delta \vx_t + \vf_{\vu_t} \delta \vu_t \right) + \bv_{t+1}\right).
\end{align*}
Using the control law~\eqref{eq:control_update_v1}, this reduces to
\begin{align*}
    &\vS_{t}\delta \vx_{t} + \bv_{t} =
    c_{\vx_t} + c_{\vx_t\vx_t} \delta \vx_t + (\bar{\blambda}_{t+1}\otimes \vf_{\vx_t\vx_t})\delta \vx_{t}
    + \vf_{\vx_t}\tran\vS_{t+1} \vf_{\vx_t}\delta \vx_t  + \vf_{\vx_{t}}\tran\bv_{t+1}\\
    &- \left( c_{\vx_t\vu_t} + \bar{\blambda}_{t+1} \otimes \vf_{\vx_t\vu_t} + \vf_{\vx_t}\tran\vS_{t+1}\vf_{\vu_t}\right)
    \left((Q_{\vu\vu}^{t})^{-1}Q_{\vu}^{t} +(Q_{\vu\vu}^{t})^{-1}Q_{\vu\vx}^{t}\delta \vx_{t} \right).
\end{align*}
Separating constant terms and terms linear in $\delta \vx_t$ yields
\begin{align*}
    \bv_t =\;& c_{\vx_t} + \vf_{\vx_{t}}\tran\bv_{t+1}
    - \left( c_{\vx_t\vu_t} + \bar{\blambda}_{t+1} \otimes \vf_{\vx_t\vu_t} + \vf_{\vx_t}\tran \vS_{t+1}\vf_{\vu_t} \right)
    \left((Q_{\vu\vu}^{t})^{-1}Q_{\vu}^{t}  \right), \\
    \vS_{t} =\;&
    c_{\vx_t\vx_t} + (\bar{\blambda}_{t+1} \otimes \vf_{\vx_t\vx_t}) + \vf_{\vx_t}\tran \vS_{t+1}\vf_{\vx_t}\\
    &- \left( c_{\vx_t\vu_t} + \bar{\blambda}_{t+1} \otimes \vf_{\vx_t\vu_t} + \vf_{\vx_t}\tran \vS_{t+1}\vf_{\vu_t} \right)
    \left((Q_{\vu\vu}^{t})^{-1}Q_{\vu\vx}^{t}  \right).
\end{align*}
Using $Q_{\vx\vx}^{t},Q_{\vx}^{t},~\text{and}~Q_{\vx\vu}^{t}$, the above equations give the backward pass:
\begin{align*}
    \bv_{t} &= Q_{\vx}^{t} - Q_{\vx\vu}^{t}(Q_{\vu\vu}^{t})^{-1}Q_{\vu}^{t}, \\
    \vS_{t} &= Q_{\vx\vx}^{t} - Q_{\vx\vu}^{t}(Q_{\vu\vu}^{t})^{-1}Q_{\vu\vx}^{t}.
\end{align*}
Finally, the terminal values of $\vS_T$ and $\bv_T$ follow from the transversality condition:
\begin{align*}
    \vS_{T} = C_{T,\vx\vx},
    \qquad
    \bv_{T} = C_{T,\vx}.
\end{align*}
This completes the derivation of the exact Newton LQR solution.
\end{proof}
It is important to note that these equations produce the same variable updates as those obtained from ``Sequential Quadratic Programming''\cite{Nocedal2006NumericalOptimization}, given that the dynamics is satisfied at the nominal: $\bar{\vx}_{t+1} = \vf(\bar{\vx}_{t},\bar{\vu}_{t})$. Thus, the resulting system of equations is equivalent to solving the optimization problem formulated in~\eqref{N-QP}. The resulting Newton algorithm is summarized in Algorithm~\ref{alg:Newton_LQR}. To improve the numerical robustness and convergence of Newton-LQR, we introduce adaptive regularization, a line-search procedure, and a modified forward pass in Algorithm~\ref{alg:Newton_LQR}.

\subsubsection{Adaptive Regularization}\label{sec:adap_reg} The control Hessian $Q_{\vu\vu}^t$, computed in line 13 of Algorithm~\ref{alg:Newton_LQR}, is regularized as

\[
Q_{\vu\vu}^t \gets Q_{\vu\vu}^t+\zeta I,
\]

where the regularization parameter is initialized at $\zeta=1$. A backward pass is accepted only if the algorithm predicts a decrease in cost. Upon acceptance, $\zeta$ is reduced by a factor of five for the subsequent iteration. Otherwise, $\zeta$ is increased by a factor of ten, and the backward pass is repeated.

\subsubsection{Line Search}
The control perturbation in line 22 is modified to incorporate a step size $\alpha$:

\[
\delta\vu_t\gets\alpha\vk_t+\vK_t\delta\vx_t.
\]

The step size starts at a value of $1$ and is progressively reduced until the ratio of the actual cost reduction to the predicted cost reduction satisfies

\[
\rho=\frac{\Delta J_{\mathrm{actual}}}
{\Delta J_{\mathrm{predicted}}}>10^{-4}.
\]

\subsubsection{Nonlinear Forward Pass} Finally, the linearized state perturbation update in line 25,

\[
\delta\hat{\vx}_{t+1}
\gets
\vf_{\vx_t}\delta\hat{\vx}_t
+
\vf_{\vu_t}\delta\vu_t,
\]

is replaced by the actual state perturbation obtained from the nonlinear dynamics:

\[
\delta\vx_{t+1}
=
\vx_{t+1}-\bar{\vx}_{t+1},
\]

where $\vx_{t+1}$ is obtained by propagating the nonlinear dynamics using the perturbed state and control. This modification avoids the approximation error introduced by linearized state propagation during the forward pass.

\begin{algorithm}
\caption{The Newton LQR solution to the Optimal Control Problem}
\label{alg:Newton_LQR}
\begin{algorithmic}[1]

\State \textbf{Initialize:} $\{\bar{\vu}_t\}_{t=0}^{T-1}$, $\{\bar{\blambda}_t\}_{t=0}^{T-1} \gets 0$,
\State Roll out the dynamics to obtain the nominal trajectory $\{\bar{\vx}_t\}_{t=0}^{T}$
\State Compute the initial cost $\cJ$

\While{not converged}
    \State Linearize the dynamics and quadratize the cost about $(\bar{\vx}_t,\bar{\vu}_t)$ according to eq.~\eqref{eq:DDP-LQR}.
        
    
    \State \textbf{Backward pass:}
    \State Initialize terminal value function derivatives:
    \[\bv_{T} \gets C_{T,\vx}, \qquad \vS_{T} \gets C_{T,\vx\vx}\]

    \For{$t = T-1$ down to $0$}
        \State
        $
        Q_{\vx}^{t} = c_{\vx_t} + \vf_{\vx_t}\tran \bv_{t+1}$
        \State $Q_{\vu}^{t} = c_{\vu_t} + \vf_{\vu_t}\tran \bv_{t+1}$
        \State $Q_{\vx\vx}^{t} = c_{\vx_t\vx_t} + \vf_{\vx_t}\tran \vS_{t+1}\vf_{\vx_t} + \bar{\blambda}_{t+1} \otimes \vf_{\vx_t\vx_t}$
        \State $Q_{\vx\vu}^{t} = c_{\vx_t\vu_t} + \vf_{\vx_t}\tran \vS_{t+1}\vf_{\vu_t} + \bar{\blambda}_{t+1} \otimes \vf_{\vx_t\vu_t}$
        \State $Q_{\vu\vu}^{t} = c_{\vu_t\vu_t} + \vf_{\vu_t}\tran \vS_{t+1}\vf_{\vu_t} + \bar{\blambda}_{t+1} \otimes \vf_{\vu_t\vu_t}$
        \State
        $
        \vk_t \gets -(Q_{\vu\vu}^{t})^{-1}Q_{\vu}^{t}
        $
        \State $
        \vK_t \gets -(Q_{\vu\vu}^{t})^{-1}Q_{\vu\vx}^{t}
        $
        \State 
        $
        \bv_t \gets Q_{\vx}^{t} - Q_{\vx\vu}^{t}(Q_{\vu\vu}^{t})^{-1}Q_{\vu}^{t}
        $
        \State $
        \vS_t \gets Q_{\vx\vx}^{t} - Q_{\vx\vu}^{t}(Q_{\vu\vu}^{t})^{-1}Q_{\vu\vx}^{t}
        $
    \EndFor

    \State \textbf{Forward pass:}
    \State Set $\vx_0 \gets \bar{\vx}_0$,  $\delta \hat{\vx}_{0} \gets \bo$
    
    \For{$t = 0$ to $T-1$}
        \State $
        \delta \vu_t \gets \vk_t + \vK_t\delta \hat{\vx}_{t}
        $
        \State $
        \vu_t \gets \bar{\vu}_t + \delta \vu_t
        $
        \State $
        \vx_{t+1} \gets \vf(\vx_t,\vu_t)
        $
        \State $
        \delta \hat{\vx}_{t+1} \gets \vf_{\vx_t}\delta \hat{\vx}_t + \vf_{\vu_t}\delta \vu_t 
        $
        \State $\bar{\blambda}_{t} = \vS_{t}\delta \hat{\vx}_{t} + \bv_{t}$
    \EndFor
    \State $\bar{\blambda}_{T} = \vS_{T}\delta \hat{\vx}_{T} + \bv_{T}$
    \State $\bar{\vx}_t \gets \vx_t$, $\bar{\vu}_t \gets \vu_t$
    

\EndWhile

\State \textbf{Return:} $\{\bar{\vx}_t\}_{t=0}^{T}, \{\bar{\vu}_t\}_{t=0}^{T-1}$

\end{algorithmic}
\end{algorithm}

\subsection{The ``modified'' Newton method - ``iLQR''}
Recall that the modified QP~\eqref{M-QP} differs
from the exact Newton QP~\eqref{N-QP} by omitting
the second-order derivatives of the equality
constraints. For the optimal control problem,
these correspond precisely to the dynamics-Hessian
terms weighted by the Lagrange multipliers.
Consequently, the modified QP gives rise to
the following LQR problem:
\begin{align}
    &\min_{\delta \vu_{t}} \sum_{t=0}^{T-1} c_{\vx_t}\delta \vx_{t} + c_{\vu_t}\delta \vu_{t} + \frac{1}{2} \delta \vx_t\tran c_{\vx_t\vx_t} \delta \vx_t + \delta \vu_t\tran c_{\vu_t\vx_t} \delta \vx_t + \frac{1}{2} \delta \vu_t\tran c_{\vu_t\vu_t} \delta \vu_t \nonumber \\
    &\quad \quad \quad \quad + C_{T,\vx}\delta \vx_{T} + \frac{1}{2} \delta \vx_{T}\tran C_{T,\vx\vx} \delta \vx_{T},\nonumber\\
    &\text{subject to:}~\delta \vx_{t+1} = \vf_{\vx_t}\delta \vx_{t} + \vf_{\vu_t}\delta \vu_{t}, \nonumber\\
    &\text{given}:~~\delta \vx_{0} = \bo.\label{eq:iLQR-LQR}
\end{align}
The solution to this LQR problems is also given by:
\begin{align*}
    \delta \vx_{t+1} =& \vf_{\vx_t}\delta \vx_{t} + \vf_{\vu_t}\delta \vu_{t},\label{eq:per_dyn_lin}\\
    \blambda_{t} =& \vS_{t}\delta \vx_{t} + \bv_{t}~~\text{with}~~\blambda_{T} = C_{T,\vx} + C_{T,\vx\vx}\delta \vx_{T},\\
    \delta \vu_{t} =& \vk_{t} + \vK_{t}\delta \vx_{t},
\end{align*}
where 
the derivatives of the augmented cost , $Q$, for this modified LQR in eq.~\eqref{eq:iLQR-LQR} are defined as follows:
\begin{align}
    Q_{\vx}^{t} =& c_{\vx_t} + \vf_{\vx_t}\tran \bv_{t+1},\\
    Q_{\vu}^{t} =& c_{\vu_t} + \vf_{\vu_t}\tran \bv_{t+1},\\
    Q_{\vx\vx}^{t} =& c_{\vx_t\vx_t} + \vf_{\vx_t}\tran \vS_{t+1} \vf_{\vx_t},\label{eq:Qxx1}\\
    Q_{\vx\vu}^{t} =& c_{\vx_t\vu_t} + \vf_{\vx_t}\tran \vS_{t+1} \vf_{\vu_t},\\
    Q_{\vu\vu}^{t} =& c_{\vu_t\vu_t} + \vf_{\vu_t}\tran \vS_{t+1} \vf_{\vu_t}.\label{eq:Quu1}
\end{align}
In contrast to the Newton LQR formulation, the modified LQR problem in eq.~\eqref{eq:iLQR-LQR} omits the tensor product terms appearing in eq.~\eqref{eq:Qxx1}-\eqref{eq:Quu1}. In particular, the terms $\bar{\blambda}_{t+1} \otimes \vf_{\vx_t\vx_t}$, $\bar{\blambda}_{t+1} \otimes \vf_{\vu_t\vx_t}$, and $\bar{\blambda}_{t+1} \otimes \vf_{\vu_t\vu_t}$ are excluded. If, in addition, the control update is computed using full perturbations in the forward pass, i.e., $\delta \vu_{t} = \vk_{t} + \vK_{t}\delta \vx_{t}$, where $\delta \vx_{t} = \vx_{t} - \bar{\vx}_{t}$, then the resulting algorithm is precisely ``iLQR''.

The absence of tensor product terms in iLQR ensures a descent direction under certain mild assumptions. Before delving into that, one should note that
the expected cost reduction for these algorithms can be simplified after substituting the control update, $\delta \vu_t = \alpha \vk_{t} + \vK_{t}\delta \vx_{t}$ and using the backward pass equations given in eq. \eqref{eq:back_pass_S} and eq.~\eqref{eq:back_pass_v} to get
\begin{equation}
    \Delta \cJ^{*} = -\left(\alpha - \alpha^{2}/2 \right)\sum_{t=0}^{T-1}(Q_{\vu}^{t})\tran (Q_{\vu\vu}^{t})^{-1}Q_{\vu}^{t},
    \label{eq:exp_cost_red}
\end{equation}
where $\alpha$ is a line search parameter used to ensure that the descent steps are within the predicted range.
One can see that $Q_{\vu\vu}^{t}$ has to be positive definite for a guarantee in expected cost reduction at each time step $t$.
The descent property of the modified QP was
established in Proposition~\ref{prop:mod_QP} under
positive-definiteness assumptions on the
quadratic objective. We now establish sufficient
conditions under which the iLQR backward pass
preserves the positive definiteness of
$Q_{\vu\vu}^{t}$ throughout the horizon,
thereby guaranteeing a descent direction.
\begin{assumption}
    $c_{\vu_t\vu_t}$ $\succ 0$. The hessian of incremental cost with respect to $\vu_{t}$ is always positive definite for all time-steps, $t$.
\end{assumption}
\begin{assumption}
    $C_{T,\vx\vx}$ $\succeq 0$. The hessian of terminal cost is positive semi-definite.
\end{assumption}
\begin{assumption}\label{ass:psd_cost}
    $\begin{bmatrix}
        c_{\vx_t\vx_t} && c_{\vx_t\vu_t}\\
        c_{\vu_t\vx_t} && c_{\vu_t\vu_t}
    \end{bmatrix}$ $\succeq 0$. The hessian matrix of incremental cost, $c(\vx_t,\vu_t)$, is positive semidefinite for all time-steps, $t$.
\end{assumption}

\begin{proposition}\label{prop:ilqr_descent}
    Under assumptions~\ref{assum:smoothness}-\ref{ass:psd_cost}, the solution to the iLQR problem is always a descent direction for the cost $c(\vx_t,\vu_t)$ and the next solution can be made feasible via a suitable line-search procedure. \\
\end{proposition}

\begin{proof}

 If we have $Q_{\vu\vu}^{t}$ to be positive definite for all time-steps then we are guaranteed a reduction in cost. At the terminal, we have
\begin{equation*}
     Q_{\vu\vu}^{T-1} = c_{\vu_{T-1}\vu_{T-1}} + \vf_{\vu_{T-1}}\tran C_{T,\vx\vx}\vf_{\vu_{T-1}}
\end{equation*}
which is positive definite following from assumptions 1 and 2. 
In addition to this, the hessian of augmented cost for iLQR is
\begin{equation}
    \begin{bmatrix}
        Q_{\vx\vx}^{T-1} & Q_{\vx\vu}^{T-1}\\
        Q_{\vu\vx}^{T-1} & Q_{\vu\vu}^{T-1}
    \end{bmatrix} = \begin{bmatrix}
        c_{\vx_{T-1}\vx_{T-1}} &  c_{\vx_{T-1}\vu_{T-1}}\\
         c_{\vu_{T-1}\vx_{T-1}} &  c_{\vu_{T-1}\vu_{T-1}} 
    \end{bmatrix} + \begin{bmatrix}
        \vf_{\vx_{T-1}}\\
        \vf_{\vu_{T-1}}
    \end{bmatrix}\tran \begin{bmatrix}
        C_{T,\vx\vx} & C_{T,\vx\vx}\\
        C_{T,\vx\vx} & C_{T,\vx\vx}
    \end{bmatrix} \begin{bmatrix}
        \vf_{\vx_{T-1}}\\
        \vf_{\vu_{T-1}}
    \end{bmatrix}.
    \label{eq:hess_matrix}
\end{equation}
In the above equation, it is evident that the matrix on left hand side is positive definite. This follows from assumptions 2 and 3, which imply:
\begin{align*}
    \begin{bmatrix}
        Q_{\vx\vx}^{T-1} & Q_{\vx\vu}^{T-1}\\
        Q_{\vu\vx}^{T-1} & Q_{\vu\vu}^{T-1}
    \end{bmatrix} &\succeq 0, \quad \text{since} \\
    \begin{bmatrix}
        c_{\vx_{T-1}\vx_{T-1}} &  c_{\vx_{T-1}\vu_{T-1}}\\
         c_{\vu_{T-1}\vx_{T-1}} &  c_{\vu_{T-1}\vu_{T-1}}
    \end{bmatrix} &\succeq 0 \quad \text{and} \quad
    \begin{bmatrix}
         C_{T,\vx\vx} & C_{T,\vx\vx}\\ 
        C_{T,\vx\vx} & C_{T,\vx\vx}
    \end{bmatrix} \succeq 0.
\end{align*}
In addition to this, $Q_{\vu\vu}^{T-1}$ is already shown to be positive definite. So, the Schur complement \cite{Gallier2019Schur} of $\begin{bmatrix}
        Q_{\vx\vx}^{T-1} & Q_{\vx\vu}^{T-1}\\
        Q_{\vu\vx}^{T-1} & Q_{\vu\vu}^{T-1}
    \end{bmatrix}$ has to be positive semidefinite, which gives
\begin{equation}
    Q_{\vx\vx}^{T-1} - Q_{\vx\vu}^{T-1}(Q_{\vu\vu}^{T-1})^{-1}Q_{\vu\vx} \succeq 0.
\end{equation}
But the above equation is nothing but the backward pass on $\vS$:
\begin{equation*}
    \vS_{T-1} =  Q_{\vx\vx}^{T-1} - Q_{\vx\vu}^{T-1}(Q_{\vu\vu}^{T-1})^{-1}Q_{\vu\vx}^{T-1} \succeq 0
\end{equation*}
Now, one can write
\begin{equation*}
    Q_{\vu\vu}^{T-2} = c_{\vu_{T-2}\vu_{T-2}} + \vf_{\vu_{T-2}}\tran \vS_{T-1}\vf_{\vu_{T-2}}
\end{equation*}
which is always positive definite as $c_{\vu_{T-2}\vu_{T-2}}$ is positive definite following from assumption 2 and $\vS_{T-1} \succeq 0$. By induction, it follows that $Q_{\vu\vu}^{t} \succ 0$ for all time-steps, $t$. This further implies that the expected cost reduction given by eq.~\eqref{eq:exp_cost_red} will always be negative guaranteeing a reduction in expected cost.
\end{proof}
Thus, this implies that the iLQR guaranties a decent direction under these mild assumptions. There are no assumptions on the dynamics, $\vf$, other than it being continuously differentiable.

In addition to this, one can ensure that the descent direction is made feasible by using the line search parameter, $\alpha$, to control the length of the step. A smaller step implies that the dynamics is closer to linear behavior and/or the perturbations are small so that the higher order terms in cost vanish.

\begin{remark}
    The convergence result in Theorem~1 of \cite{Wang2025SearchFeedbackRL} extends to the present setting, implying that iLQR converges to a stationary point when equipped with a sufficiently small line search parameter, $\alpha$.
\end{remark}

\subsection{The ``approximate'' Newton methods - ``DDP'' and ``Stagewise-Newton''}
Thus far, we have focused on the Newton QP in eq.~\eqref{N-QP} and the modified QP in eq.~\eqref{M-QP} as formulations for solving the optimal control problem defined by eq.~\eqref{eq:ocp} and eq.~\eqref{eq:dynamics}. We now consider other related algorithms that can be interpreted as approximations to the Newton QP for the solution of the same optimal control problem.

An alternative interpretation of the Newton LQR problem in eq.~\eqref{eq:DDP-LQR} leads to the equations commonly associated with DDP. Recall that the costate is represented in the affine form
\[
\blambda_{t} = \bv_{t} + \vS_{t}\delta \vx_{t}.
\]
At the same time, the costate may also be expressed as
\[
\blambda_{t} = \bar{\blambda}_{t} + \delta \blambda_{t},
\]
where $\bar{\blambda}_{t}$ denotes the costate from the previous iteration. This suggests the identification
\[
\bar{\blambda}_{t} = \bv_{t}, \qquad \delta \blambda_{t} = \vS_{t}\delta \vx_{t},
\]
so that
\begin{align*}
    \blambda_{t} = \bar{\blambda}_{t} + \delta \blambda_{t} = \bv_{t} + \vS_{t}\delta \vx_{t}.
\end{align*}
This interpretation leads directly to the DDP equations when substituted into the Newton LQR solution. Accordingly, the distinction between the backward pass of DDP and that of Newton LQR lies in the following three terms:
\begin{align*}
    Q_{\vx\vx}^{t} =& c_{\vx_t\vx_t} + \vf_{\vx_t}\tran \vS_{t+1}\vf_{\vx_t} + \bv_{t+1} \otimes \vf_{\vx_t\vx_t};\\
    Q_{\vx\vu}^{t} =& c_{\vx_t\vu_t} + \vf_{\vx_t}\tran \vS_{t+1}\vf_{\vu_t} + \bv_{t+1} \otimes \vf_{\vx_t\vu_t};\\
    Q_{\vu\vu}^{t} =& c_{\vu_t\vu_t} + \vf_{\vu_t}\tran \vS_{t+1}\vf_{\vu_t} + \bv_{t+1} \otimes \vf_{\vu_t\vu_t}.
\end{align*}
The principal distinction between Newton-LQR
and DDP lies in the multipliers used to weight
the dynamics-Hessian terms. Newton-LQR employs
the Lagrange multipliers inherited from the
SQP iteration, whereas DDP uses the gradient
of its locally approximated value function.
Consequently, DDP replaces
$\bar{\blambda}_{t+1}$ with $\bv_{t+1}$
in the second-order dynamics terms.
It is important to note that, unlike Newton LQR, DDP does not require the computation of $\bar{\blambda}_{t}$ through forward pass, since this quantity is not treated separately and is same as $\bv_t$ which is already computed in the backward pass.\\

\begin{remark}
\label{rem:ddp_newton_consistency}
Under Assumption~\ref{assum:smoothness}, at a feasible stationary trajectory $(\vx^*,\vu^*)$ satisfying the first-order necessary conditions for optimality, the DDP feedforward update vanishes, i.e., $\mathbf{k}_t=\boldsymbol{0}$. Consequently, the DDP value-function gradient recursion coincides with the costate recursion, yielding $\bv_t=\boldsymbol{\lambda}_t^*$ for all $t=0,\ldots,T$. Thus, the dynamics-Hessian terms in DDP and the exact Newton/SQP Hessian coincide at stationarity. Away from stationarity, this equivalence is not guaranteed as will be seen empirically in Section 5.
\end{remark}

Next, we move on to derive the ``Stagewise-Newton'' equations. To do so, first consider the following equation from Newton LQR:
\begin{align*}
    \bar{\blambda}_{t}^{\text{iter}} = \underbrace{\vS_{t}\delta \hat{\vx}_{t} + \bv_{t},}_{\text{at iter-1}}
\end{align*}
where iter denotes value at the current iteration and (iter-1) denotes values at the previous iteration.
Using backward pass from eq.~\eqref{eq:back_pass_S} and eq.~\eqref{eq:back_pass_v}, we have
\begin{align*}
    \bar{\blambda}_{t}^{\text{iter}} = \underbrace{\left( Q_{\vx\vx}^{t} - Q_{\vx\vu}^{t}(Q_{\vu\vu}^{t})^{-1}Q_{\vu\vx}^{t}\right)\delta \hat{\vx}_t + Q_{\vx}^{t} - Q_{\vx\vu}^{t}(Q_{\vu\vu}^{t})^{-1}Q_{\vu}^{t}}_{\text{all derivatives at iter-1}},
\end{align*}
which can be written as:
\begin{align*}
    \bar{\blambda}_{t}^{\text{iter}}  = \underbrace{Q_{\vx}^{t} + Q_{\vx\vx}^{t}\delta \hat{\vx}_{t} + Q_{\vx\vu}^{t}\left( (-Q_{\vu\vu}^{t})^{-1}Q_{\vu\vx}^{t}\delta \hat{\vx}_{t} + (-Q_{\vu\vu}^{t})^{-1}Q_{\vu}^{t} \right)}_{\text{all derivatives at iter-1}},
\end{align*}
which further reduces to
\begin{align*}
    \bar{\blambda}_{t}^{\text{iter}}  =& c_{\vx_{t}} + \vf_{\vx_t}\tran\bv_{t+1} + c_{\vx_{t}\vx_t}\delta \hat{\vx}_{t} + \vf_{\vx_t}\tran\vS_{t+1}\vf_{\vx_t}\delta \hat{\vx}_t\\
    &+  \underbrace{(\bar{\blambda}_{t+1}\otimes \vf_{\vx_t\vx_t}) \delta \hat{\vx}_{t} + (c_{\vx_t\vu_t} + \vf_{\vx_t}\tran\vS_{t+1}\vf_{\vu_t} + (\bar{\blambda}_{t+1}\otimes \vf_{\vx_t\vu_t}))\delta\vu_{t}}_{\text{all derivatives at iter-1}}.
\end{align*}
Upon rearranging terms, we get
\begin{align}
    \bar{\blambda}_{t}^{\text{iter}}  =& c_{\vx_{t}} + c_{\vx_{t}\vx_t}\delta \hat{\vx}_{t} + c_{\vx_{t}\vu_t}\delta \vu_{t}  + \vf_{\vx_t}\tran\left( \vS_{t+1} \left( \vf_{\vx_t}\delta \hat{\vx}_t + \vf_{\vu_t}\delta \vu_{t}\right) + \bv_{t+1}\right) \nonumber\\
    &+  \underbrace{(\bar{\blambda}_{t+1}\otimes \vf_{\vx_t\vx_t}) \delta \hat{\vx}_{t} + (\bar{\blambda}_{t+1}\otimes \vf_{\vx_t\vu_t})\delta\vu_{t}}_{\text{all derivatives at iter-1}},\label{eq:lin_SN}
\end{align}
which under linearity assumptions can be written as:
\begin{align}
    \bar{\blambda}_{t}^{\text{iter}} = \underbrace{c_{\vx_t} + \vf_{\vx_t}\tran \bar{\blambda}_{t+1}}_{\text{all values at iter}}
    \label{eq:SN_back_pass1}
\end{align}
If rather than computing $\bar{\blambda}_{t}$ in the forward pass as in Newton LQR, we compute this in backward pass using eq.~\eqref{eq:SN_back_pass1}, we will get the ``Stagewise-Newton''. Rest of the equations remain the same for Stagewise-Newton. One should note here that the expression in eq.~\eqref{eq:SN_back_pass1} is a necessary condition of the original optimal control problem and should be valid only at the optimal. This equality holds in an approximate sense away from the optimal.\\
Unlike DDP and iLQR, Stagewise-Newton also uses the linearized perturbation dynamics in the forward pass to update the control. It is important to note that the linearity assumption employed in deriving eq.~\eqref{eq:SN_back_pass1} is often a poor approximation and may result in slow convergence. Specifically, the $\bar{\blambda}_{t+1}$ terms in the tensor product in eq.~\eqref{eq:lin_SN} are evaluated using information from the previous iteration (iter-1), whereas for the linearity assumption to be valid, these terms should be computed at the current iteration (iter). Thus, the transition from eq.~\eqref{eq:lin_SN} to eq.~\eqref{eq:SN_back_pass1} relies on an approximation whose accuracy is limited. In addition to this, taking linear perturbation dynamics leads to slow convergence as will be discussed in Section~\ref{sec:algorithms}.
\begin{remark}
    The preceding analysis indicates that neither DDP nor Stagewise-Newton constitutes a true Newton method. Rather, both are obtained by introducing approximations to the true Newton LQR problem.
\end{remark}

\subsection{The Hybrid Newton Approach}
\label{subsec:hybrid_newton}

Under appropriate regularity and second-order sufficient
conditions, Newton-type methods and DDP exhibit locally
quadratic convergence. To characterize this convergence
rate, let
$\vU^{(k)}=\{\vu_t^{(k)}\}_{t=0}^{T-1}$
denote the control trajectory at iteration $k$, and let
$\vU^*=\{\vu_t^*\}_{t=0}^{T-1}$
denote a locally optimal control trajectory. Quadratic
convergence is characterized by
\begin{align}
    \left\|
    \vU^{(k+1)}-\vU^*
    \right\|_2
    \leq
    c
    \left\|
    \vU^{(k)}-\vU^*
    \right\|_2^2,
    \label{eq:con_map}
\end{align}
where $c>0$ is a constant independent of $k$, and
the inequality holds when the iterates are sufficiently
close to $\vU^*$.

Although the inclusion of dynamics Hessian information
enables Newton-type methods and DDP to achieve faster
local convergence, it does not necessarily improve their
performance when the iterates are far from a stationary
solution. As established in the preceding sections,
the dynamics Hessian terms can render the local quadratic
subproblem indefinite, necessitating regularization
to obtain a descent direction. In contrast, iLQR
omits these terms and, under suitable convexity
assumptions on the cost functions, guarantees a descent
direction without requiring regularization.

These observations motivate a hybrid optimization
strategy that exploits the complementary advantages
of iLQR and second-order methods. The proposed
approach consists of two phases. During the initial
phase, iLQR is employed to exploit its favorable
descent properties and obtain a trajectory in the
vicinity of a stationary solution. Once sufficient
progress has been achieved, the algorithm transitions
to either DDP or the proposed Newton-LQR method
to exploit the additional second-order dynamics
information and potentially accelerate convergence.

To determine when the transition should occur, we
consider the relative change in the objective function
between successive iterations. Specifically, the
switching criterion is defined as
\begin{align}
    \frac{
    \left|\cJ_{k+1}-\cJ_k\right|
    }{
    \max\left\{\left|\cJ_k\right|,\varepsilon\right\}
    }
    <\psi,
    \label{eq:hybrid_switching}
\end{align}
where $\cJ_k$ denotes the objective function value
at iteration $k$, $\psi>0$ is a user-defined threshold,
and $\varepsilon>0$ is a small constant introduced
to prevent division by zero.

Once the criterion in
\eqref{eq:hybrid_switching} is satisfied, the
iLQR backward pass is replaced by the corresponding
DDP or Newton-LQR backward pass, thereby incorporating
the dynamics Hessian terms into the local quadratic
model. The second-order method is retained for the
remainder of the optimization.

It is important to note that a small relative change
in the objective function does not necessarily imply
proximity to a stationary solution, as it may also
indicate stagnation. Therefore,
\eqref{eq:hybrid_switching} is employed as a
practical switching heuristic rather than a
theoretical certificate of stationarity.

The objective of this hybrid strategy is to combine
the reliable early-stage convergence behavior of
iLQR with the potentially faster local convergence
of second-order methods. In the numerical section,
we investigate the stage of optimization at which
activating the dynamics Hessian information becomes
advantageous.

\section{Algorithmic Performance: Comparison on Benchmark Examples}\label{sec:algorithms}
Having examined Newton's method and its connections to iLQR, DDP, and Stagewise-Newton through the lens of the Newton QP formulation, we now turn to a numerical comparison of these methods to evaluate their practical performance and convergence behavior. In this section, we present a direct comparison of these algorithms on representative control tasks. By examining convergence rates, regularization, and step sizes, we aim to expose the practical advantages and limitations that arise solely from their respective local modeling choices.  

\subsection{Benchmark Tasks for Algorithmic Comparison}
Six problems are considered for this comparison. For the first four problems, both the running and terminal costs are quadratic. The general cost form is given by:
\begin{align*}
    J = \min_{\vu(t)} \frac{1}{2}\int_{0}^{t_f} (\vx-\vx_f)\tran \vQ(\vx - \vx_f) + \vu\tran \vR\vu ~~~dt \\
    + \frac{1}{2} (\vx_{t_f} - \vx_f)\tran \vQ_f (\vx_{t_f} - \vx_f).
\end{align*}
The parameters for each of these problems are summarized in Table~\ref{tab:sys_params}. The governing equations were discretized using a fourth-order Runge–Kutta (RK4) integration scheme. In the table, $\Delta t$ denotes the duration of the zero-order hold used to solve the problem. The cost specifications for the last two problems are provided in their respective descriptions.

\subsubsection{Pendulum Swing-up}
The first problem considered is the pendulum swing-up task. The pendulum dynamics is governed by the following:
\begin{align*}
\frac{d}{dt}\begin{bmatrix}
    \theta\\
    \dot{\theta}
\end{bmatrix}
=
\begin{bmatrix}
    \dot{\theta} \\
    -\frac{g\sin(\theta)}{L} + \frac{u}{mL^{2}}
\end{bmatrix}
\end{align*}
where $\theta$ and $\dot{\theta}$ are the states, and $u$ denotes the control moment applied at the other end of the rod. $m$ and $L$ are the mass of the bob and the length of the massless rod, respectively, and $g$ denotes gravity.  Table~\ref{tab:sys_params} summarizes the parameters used for discretization of the dynamics. 

\subsubsection{Cart-pole Swing-up}
The dynamics governing this problem are:
\begin{align*}
    \frac{d}{dt}\begin{bmatrix}
        x\\
        \dot{x}\\
        \theta\\
        \dot{\theta}
    \end{bmatrix}~=~\begin{bmatrix}
        \dot{x}\\
        \frac{u + ml\dot{\theta}^{2}\sin(\theta) - mg\sin(\theta)\cos(\theta)}{M + m\sin^{2}(\theta)}\\
        \dot{\theta}\\
        \frac{-u\cos(\theta) - ml\dot{\theta}^{2} \sin(\theta)\cos(\theta) + (M + m)g\sin(\theta)}{l(M + m\sin^{2}(\theta))}
    \end{bmatrix},
\end{align*}
where $x$ and $\dot{x}$ denote the position and velocity of the cart, $\theta$ and $\dot{\theta}$ denote the angle and angular velocity of the pole, and $u$ is the control input corresponding to the force applied to the cart. $M$, $m$ and $l$ are the masses of cart, bob, and length of the pole, respectively. As in the pendulum task, the cost and discretization parameters are summarized in Table~\ref{tab:sys_params}.

\subsubsection{Test Problem-1}
A hypersensitive problem \cite{Rao2000} with dynamics given by
\begin{equation*}
    \dot{\vx} = -\vx^{3} + \vu
\end{equation*}
is considered and referred to as ``Test Problem-1''. 

\subsubsection{Test Problem-2}
Another one-dimensional test problem is taken from \cite{jacobson1970ddp}. The governing equation for this problem is
\begin{equation*}
    \dot{\vx} = -0.2\vx + 10\tanh(\vu).
\end{equation*}
This problem is labeled ``Test Problem-2''. 
\begin{table*}[!htpb]
    \centering
    \tiny
    \begin{tabular}{|c|c|c|c|c|c|c|c|c|}
        \hline
        System & $ t_f$(s) & $\Delta t$(s) & T & $\vQ$ & $\vR$ & $\vQ_f$ & $\vx_{0}$ & $\vx_{f}$\\
        \hline
        Pendulum    & 5      &  0.1 & 50 & 3$I_{2\times 2}$ &  3 & 30$I_{2\times 2}$  & $[0, 0]'$ & $[\pi,0]'$     \\
        \hline
        Cartpole    & 3       &  0.1 & 30 & 10$I_{4\times 4}$ & 0.001 & 10000$I_{4\times 4}$ &  $[0,0,\pi,0]'$ & $[0,0,0,0]'$  \\
        \hline
        TP-1    & 25       &  0.01 & 2500 & 2 & 2 & $10^{10}$ & 1 & 1.5    \\
        \hline
        TP-2 & 0.5 &  $10^{-4}$ & 5000 & 20 & 2 & 20 & 5 & 0\\
        \hline
    \end{tabular}
    \caption{Parameters for all four problems. TP-1 denotes Test problem-1 and TP-2 denotes Test problem-2.}
    \label{tab:sys_params}
\end{table*}

\subsubsection{Test Problem - 3}
In addition to these problems, we also consider two problems from \cite{LiaoShoemaker1992}. The first of these is labeled ``Test Problem-3'' and is described by:
\begin{align*}
    \min \mathcal{J} =& \sum_{t=1}^{T-1}\left( \sum_{i=1}^{n} \left( (\vx_{t})_{i} + \frac{1}{4} \right)^{4} + \sum_{j=1}^{m} \left( (\vu_{t})_{j} + \frac{1}{4} \right)^{4} \right) + \sum_{i=1}^{n} \left( (\vx_{T})_{i} + \frac{1}{4} \right)^{4},\\
    &\text{where}~~ \vx_{t+1} = \mathbf{A}\vx_{t} + \mathbf{B}\vu_{t} + (\vx_{t}\tran \mathbf{C} \vu_{t}) \boldsymbol{\gamma} \quad t = 1,2,\cdots T-1\\
    &\vx_{1} = \boldsymbol{0}_{n\times 1} \in \mathbb{R}^{n}\\
    &\mathbf{A} \in \mathbb{R}^{n \times n}~~\text{and}~~(\mathbf{A})_{ij} = \begin{cases}
        0.5 \quad \text{if}~i=j\\
        0.25 \quad \text{if}~j=i+1\\
        -0.25 \quad \text{if}~j=i-1\\
        0\quad \text{elsewhere}
    \end{cases}~~i=1,\cdots n\\
    & \mathbf{B} \in \mathbb{R}^{n \times m}, ~~(\mathbf{B})_{ij} = \frac{i-j}{m+n},~~i=1,\cdots n ~~\text{and}~~j= 1,\cdots m\\
    & \mathbf{C} \in \mathbb{R}^{n \times m}, ~~(\mathbf{C})_{ij} = \mu\frac{(i+j)}{m+n},~\mu~\text{is a scalar},\\
    & \boldsymbol{\gamma} \in \mathbb{R}^{n}~\text{ and }~ \boldsymbol{\gamma} = \boldsymbol{1}_{n \times 1}.
\end{align*}
In the above problem, $(\cdot)_{ij}$ denotes the element of a matrix at the $i$-th row and $j$-th column.

\subsubsection{Test Problem-4}
The sixth and final problem is labeled ``Test Problem-4'' and is described below\cite{LiaoShoemaker1992}:
\begin{align*}
    \min \mathcal{J} =& \sum_{t=1}^{T-1} || \vx_{t} ||_{2}^{2} \left( \sin \left( \frac{||\vu_{t}||_{2}^{2}}{m} \right) + 1 \right)  + || \vx_{T}||_{2}^{2}\\
    &\text{where}~(\mathbf{x}_{t+1})_{i} = \sin \left( (\mathbf{x}_{t})_{i} \right) + \left( \mathbf{F} \mathbf{W}(\mathbf{u}_{t}) \right)_{i}\\
    &(\cdot)_{i}~\text{denotes the }~i\text{-th component of the vector}~(\cdot)\\
    &(\mathbf{x}_{1})_{i} = \frac{i}{2n}\quad i = 1,\cdots n\\
    &(\mathbf{F})_{ij} = \frac{i+j}{2n}\quad i = 1, \cdots n~~\text{and}~~j=1,\cdots m\\
    &\mathbf{W}(\mathbf{u}_{t}) = \left( \sin\left(\left(\mathbf{u}_{t}\right)_{1}\right),\cdots,\sin\left(\left(\mathbf{u}_{t}\right)_{m}\right)\right)\tran \in \mathbb{R}^{m}\\
    &\vx_{t} \in \mathbb{R}^{n} \quad \text{and} \quad \vu_{t} \in \mathbb{R}^{m}.
\end{align*}


\subsection{The impact of Forward-Pass}
A critical distinction between the Stagewise-Newton method and the DDP/iLQR family of algorithms lies in the forward pass, specifically in how state perturbations are propagated to compute the control update. In DDP and iLQR, the feedback term depends on the full nonlinear state perturbation, $\delta \vx_t$. The control update is
$
\delta \vu_t = \alpha \vk_t + \vK_{t} \delta \vx_t
$,
where $\delta \vx_t$ is the difference between the new and previous (nominal) trajectories:
\begin{align*}
    \delta \vx_t = \vf(\vx_{t}, \vu_{t}) - \vf(\bar{\vx}_{t}, \bar{\vu}_{t}).
\end{align*}
In contrast, the Stagewise-Newton method employs a control update, $\delta \vu_t = \alpha \vk_t + \vK_{t} \delta \hat{\vx}_t$, that relies on a linearized propagation of the state perturbation, $\delta \hat{\vx}_t$:
\begin{align*}
    \delta \hat{\vx}_{t+1} = \vf_{\vx_t}\delta \hat{\vx}_{t} + \vf_{\vu_t} \delta \vu_{t}.
\end{align*}
The performance implications of this choice were highlighted by Liao and Shoemaker \cite{LiaoShoemaker1992}, who observed that the linearized forward pass of the Stagewise-Newton method leads to significantly slower convergence compared to DDP. To isolate the effect of the forward pass, they introduced a ``Mixed" algorithm combining a DDP-style forward pass with a Stagewise-Newton backward pass, confirming that using the full nonlinear state perturbation, $\delta \vx_t$, accelerates convergence.

To analyze this difference quantitatively, we examine the Taylor series expansion of the state perturbation. Since the initial state is fixed, $\delta \vx_0 = \boldsymbol{0}$, the first control perturbation is purely feedforward: $\delta \vu_0 = \alpha \vk_0$. The resulting state perturbation at the next step, $\delta \vx_1$, can be expanded as:
\begin{dmath*}
    \delta \vx_{1} = \vf_{\vx_{0}}\delta \vx_{0} + \vf_{\vu_{0}}\delta \vu_{0} + \frac{1}{2} \delta \vx_{0}\tran \vf_{\vx_{0}\vx_{0}} \delta \vx_{0}+ \delta \vu_{0}\tran \vf_{\vu_{0}\vx_{0}} \delta \vx_{0}+ \frac{1}{2} \delta \vu_{0}\tran \vf_{\vu_{0}\vu_{0}} \delta \vu_{0} + \alpha^{3}(\text{third order terms}) + \alpha^{4}(\text{fourth order terms}) + \cdots \cdots
\end{dmath*}
Substituting the initial conditions ($\delta \vx_{0} = \bo$) and the control update ($\delta \vu_{0} = \alpha \vk_{0} + \vK_{0}\delta \vx_{0}$) yields an expansion in powers of the line search parameter $\alpha$:
\begin{dmath*}
    \delta \vx_{1} = \alpha \vf_{\vu_{0}} \vk_{0} + \frac{1}{2} \alpha^{2} \vk_{0}\tran \vf_{\vu_{0}\vu_{0}} \vk_{0} + + \alpha^{3}(\text{third order terms}) + \alpha^{4}(\text{fourth order terms}) + \cdots \cdots
\end{dmath*}
This recursive structure propagates through time. The state perturbation $\delta \vx_t$ at any time step $t$ is an expansion in powers of $\alpha$, where the coefficients depend on the system Jacobians and the feedforward control modifications $\{\vk_i\}_{i=0}^{t-1}$. The general form of the first-order term in $\alpha$ is:
\begin{dmath}
    \delta \vx_{t} = \alpha \left[ \vf_{\vu_{t-1}}\vk_{t-1} + \sum_{j=0}^{t-2} \left( \prod_{i=j+1}^{t-1} (\vf_{\vx_{i}} + \vf_{\vu_{i}}\vK_{i}) \right) \vf_{\vu_{j}}\vk_{j} \right] +  \alpha^{2}(\text{second order terms}) + \alpha^{3}(\text{third order terms}) + \cdots \cdots \label{eq:lin_exp_dyn}
\end{dmath}

\begin{remark}
In the perturbation expansion of the dynamics, each $k$-th order term is scaled by a factor of $\alpha^k$. Consequently, one can diminish the influence of higher-order terms and recover the linearized dynamics by reducing $\alpha$.
\end{remark}

The Stagewise-Newton method, by design, retains only the first-order terms in $\alpha$ from this expansion for its feedback correction. This corresponds to an implicit assumption that higher-order dynamics are negligible, which can be a poor approximation. In contrast, DDP, iLQR, and Newton-LQR use a line search on $\alpha$ to actively manage the influence of these higher-order terms, reducing $\alpha$ only when necessary to ensure a sufficient decrease in the cost function. This adaptive approach avoids prematurely discarding valuable nonlinear information.

\begin{figure}[!htbp]
    \centering
    \subfloat[$\alpha=1$]{\includegraphics[width=0.48\linewidth]{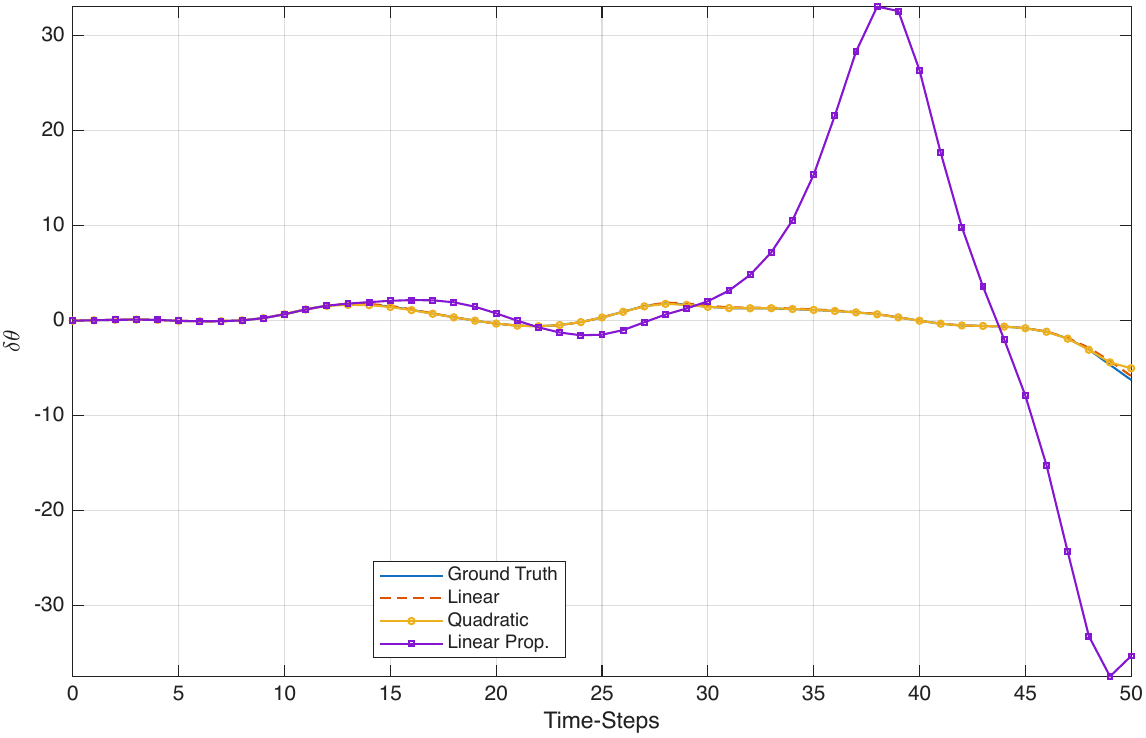}\label{subfig:alpha1}}
   \hfill
   \subfloat[$\alpha=0.1$]{\includegraphics[width=0.48\linewidth]{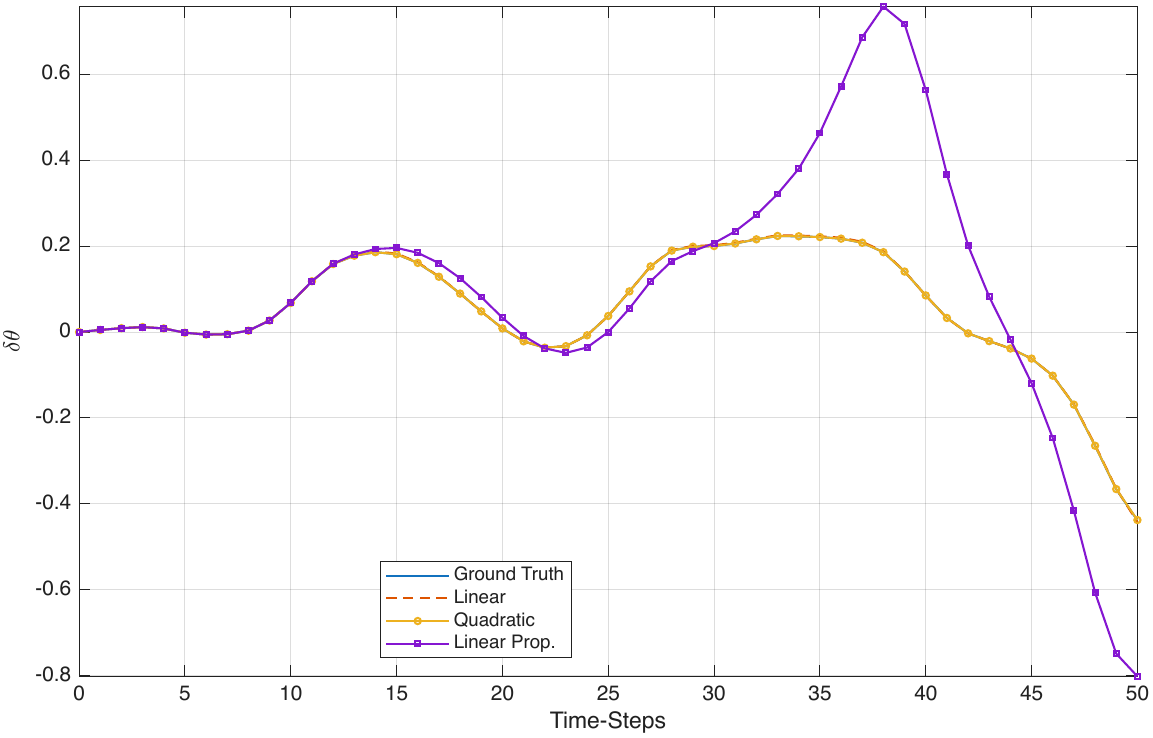}\label{subfig:alpha01}}%
   \vfill
   \subfloat[$\alpha=0.05$]{\includegraphics[width=0.48\linewidth]{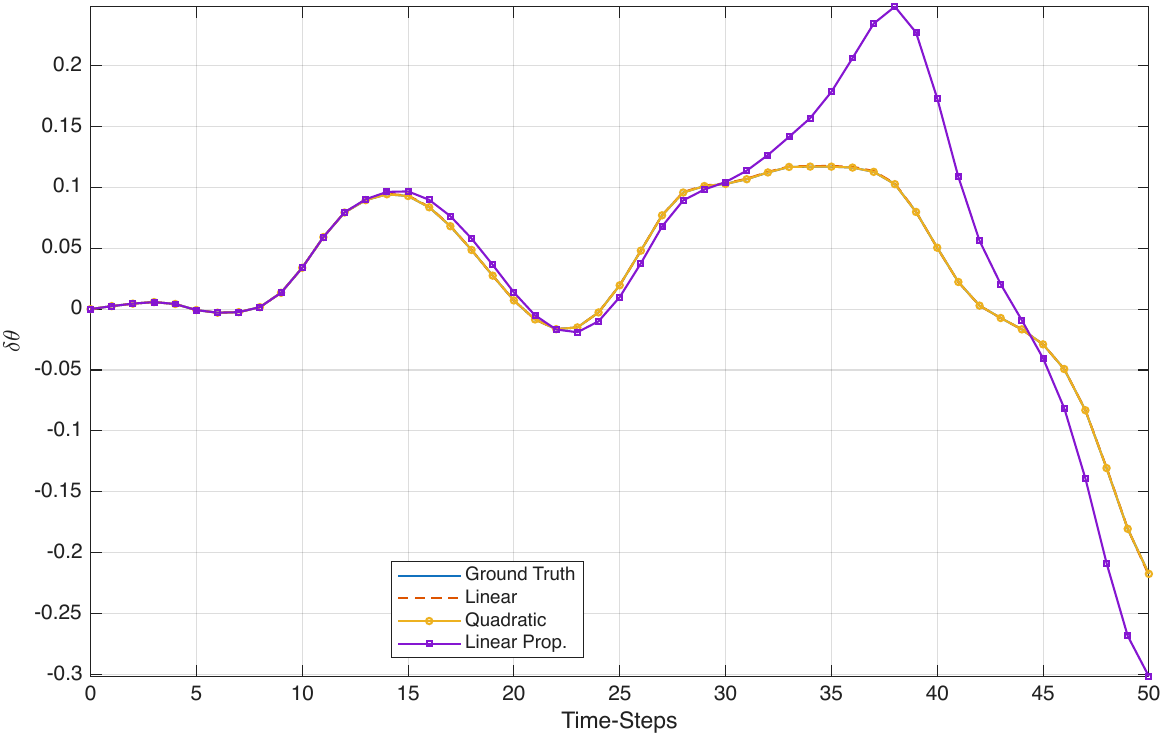}\label{subfig:alpha005}}
   \hfill
   \subfloat[$\alpha=0.01$]{\includegraphics[width=0.48\linewidth]{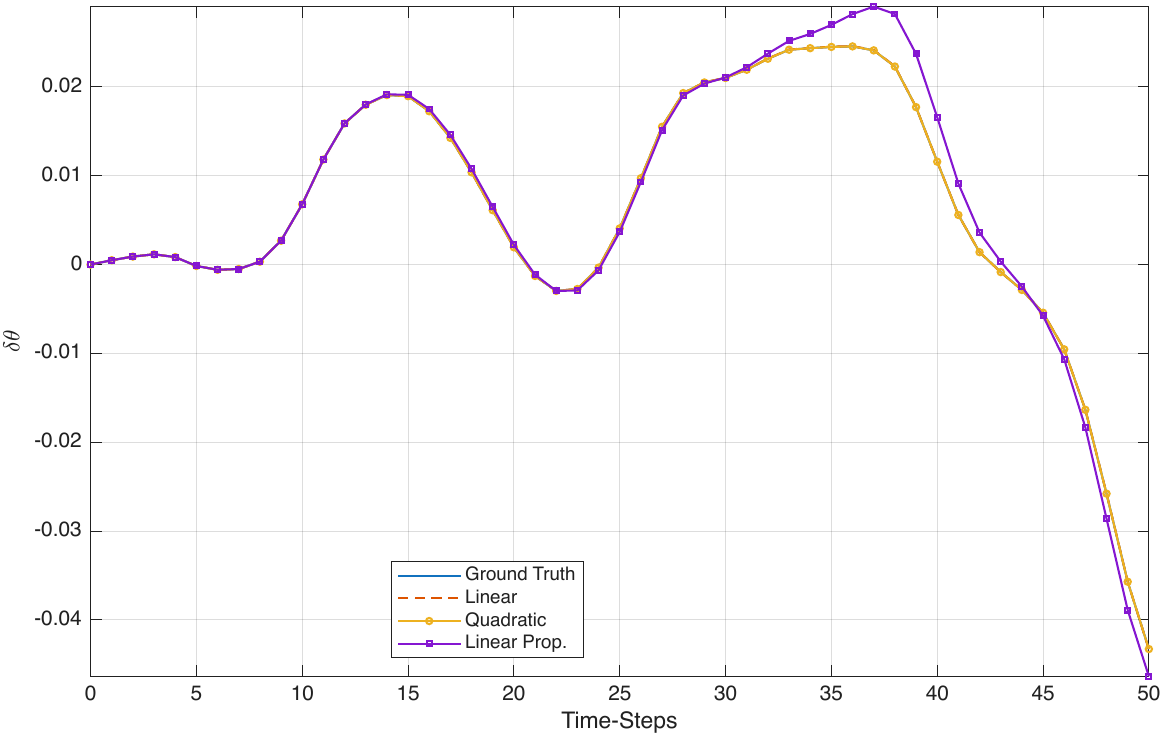}\label{subfig:alpha001}}
    \caption{Propagation of the state perturbation, $\delta\theta$, for the pendulum task. The ground truth (solid line) is compared with predictions from a local linear model (dashed), a local quadratic model (circles), and the recursive Stagewise-Newton update (squares). The Stagewise-Newton update only converges to the ground truth for very small line search parameters ($\alpha \ll 1$), whereas the local models remain accurate, highlighting the importance of using the full state perturbation in the forward pass.}
    \label{fig:alpha_comp_SN}
\end{figure}

This theoretical difference is demonstrated empirically in Fig.~\ref{fig:alpha_comp_SN} for the pendulum swing-up task. The figure compares four methods of propagating the state perturbation $\delta\theta$: 
1) the ground truth (full nonlinear difference), 2) a local linear approximation, 3) a local quadratic approximation, and 4) the recursive linear propagation used in Stagewise-Newton. For a given iteration, the local linear and quadratic models are nearly indistinguishable from the ground truth, indicating the dynamics are well-behaved locally. However, the Stagewise-Newton propagation (``Linear Prop.") significantly deviates from the ground truth, especially for larger $\alpha$. This error reduces as $\alpha$ is reduced (e.g., $\alpha=0.01$), confirming that the method's validity is restricted to small step sizes.

These findings on the forward pass directly correlate with convergence speed, as shown in Table~\ref{tab:iterations_to_optimal}, where we tested these algorithms for ``Test Problem-4''. The data, including results for DDP and Stagewise-Newton from \cite{LiaoShoemaker1992}, demonstrates that iLQR, DDP and Newton-LQR significantly outperform the Stagewise-Newton method due to their handling of nonlinearity. 

\begin{table}[!htbp]
    \centering
    \caption{Number of iterations to convergence for Test Problem 4. $T=100$, $n=100$, and $m=10$.}
    \label{tab:iterations_to_optimal}
    \begin{tabular}{|l|c|c|c|c|}
        \hline
        & \multicolumn{4}{c|}{Number of Iterations} \\
        \cline{2-5}
        Starting Point\textsuperscript{a} & iLQR & DDP\textsuperscript{b} & Stagewise-Newton\textsuperscript{b} & Newton-LQR  \\
        \hline
        1 & 4 & 8 & 132 & 4\\
        2 & 4 & 9 & 127 & 6\\
        3 & 5 & 10 & 128 & 6\\
        4 & 4 & 10 & 134 & 4\\
        5 & 5 & 10 & 134 & 4\\
        \hline
        \multicolumn{4}{l}{\textsuperscript{a}\footnotesize{Initial control guesses are defined in Table~\ref{tab:starting_points}.}} \\
        \multicolumn{4}{l}{\textsuperscript{b}\footnotesize{Data for DDP and Newton methods from \cite{LiaoShoemaker1992}.}}
    \end{tabular}
\end{table}

\begin{table}[htbp]
    \centering
    \caption{Initial control trajectories for Test Problem 4 ($T=100, m=10$).}
    \label{tab:starting_points}
    \begin{tabular}{|c|l|}
        \hline
        \textbf{Starting Point} & \textbf{Initial Control}, $\vu_t$ \\
        \hline
        1 & $\vu_{t} = \boldsymbol{0}$\\
        \hline
        2 & $\vu_{t} = 0.01 \cdot \boldsymbol{1} $\\
        \hline
        3 & $\vu_t = -0.01 \cdot \boldsymbol{1} $ \\
        \hline
        4 &
        $\begin{cases}
          \vu_{t} = 0.01 \cdot \boldsymbol{1} & t \text{ is odd} \\
          \vu_{t} = \boldsymbol{0}  & t \text{ is even}
        \end{cases}
        $ \\
        \hline
        5 &
        $\begin{cases}
          \vu_{t} = -0.01 \cdot \boldsymbol{1} & t \text{ is odd} \\
          \vu_{t} = \boldsymbol{0}   & t \text{ is even}
        \end{cases}
        $ \\
        \hline
    \end{tabular}
\end{table}

\subsection{The Need for Regularization in the Presence of Dynamics Hessians}
The preceding discussion demonstrates the advantages of using full state perturbations rather than the linearized approximations employed in Stagewise-Newton. We therefore adopt this approach in the subsequent numerical investigations, restricting our attention to iLQR, DDP, and Newton-LQR.

In this subsection, we examine the descent directions predicted by iLQR, and contrast it with DDP and Newton-LQR. To begin with, we present the expected reduction in cost predicted at different time steps of iLQR and DDP. The net ``expected cost reduction" in these algorithms is given by eq.~\eqref{eq:exp_cost_red}. It should be noted that the term $Q_{\vu\vu}^{t}$ should be positive definite to guarantee a certain reduction in cost. From Proposition \ref{prop:ilqr_descent}, it can stated that under some mild assumptions, this term is always positive definite for iLQR as described in the previous section. However, it may become negative-definite/indefinite for DDP and Newton-LQR. If we do a similar analysis for the cost hessian like iLQR ( (eq.~\eqref{eq:hess_matrix} in Proposition \ref{prop:ilqr_descent}) for DDP and Newton-LQR, we get
\begin{align*}
    \begin{bmatrix}
        Q_{\vx\vx}^{T-1} & Q_{\vx\vu}^{T-1}\\
        Q_{\vu\vx}^{T-1} & Q_{\vu\vu}^{T-1}
    \end{bmatrix} = 
    \begin{bmatrix}
        c_{\vx_{T-1}\vx_{T-1}} &  c_{\vx_{T-1}\vu_{T-1}}\\
         c_{\vu_{T-1}\vx_{T-1}} &  c_{\vu_{T-1}\vu_{T-1}} 
    \end{bmatrix} +& \begin{bmatrix}
        \vf_{\vx_{T-1}}\\
        \vf_{\vu_{T-1}}
    \end{bmatrix}\tran \begin{bmatrix}
        C_{T,\vx\vx} & C_{T,\vx\vx}\\
        C_{T,\vx\vx} & C_{T,\vx\vx}
    \end{bmatrix} \begin{bmatrix}
        \vf_{\vx_{T-1}}\\
        \vf_{\vu_{T-1}}
    \end{bmatrix}\\
    &+ \blambda_{T} \otimes 
    \begin{bmatrix}
        \vf_{\vx_{T-1}\vx_{T-1}} & \vf_{\vx_{T-1}\vu_{T-1}}\\
        \vf_{\vu_{T-1}\vx_{T-1}} & \vf_{\vu_{T-1}\vu_{T-1}}
    \end{bmatrix}.
\end{align*}
\begin{remark}
    The positive-definiteness guarantee established for iLQR does not
    generally extend to DDP and Newton-LQR due to the presence of
    dynamics Hessian terms. Consequently,
    regularization may be required to restore positive definiteness
    and obtain a descent direction for the local quadratic model.
\end{remark}
The inclusion of dynamics Hessian terms in DDP and Newton-LQR
can render $Q_{\vu\vu}^t$ indefinite, necessitating regularization
to obtain a descent direction. This analysis is presented in~\ref{app:regularization},
where we examine the effect of regularization on the predicted
cost reduction and show that for sufficiently large regularization, a descent direction can always be obtained, and under heavy regularization,
the control update asymptotically approaches a scaled
gradient-descent step.

\begin{table}[!htbp]
    \centering
    \caption{Positive definiteness of $Q_{\vu\vu}^t$ across the benchmark problems. $\gamma$ denotes the minimum eigenvalue of $Q_{\vu\vu}^{t}$ over all time, $t$. $N$ is the number of times $Q_{\vu\vu}^{t}$ is non positive-definite over the horizon, $T$.}
    \label{tab:Quu_PD}
    \begin{tabular}{|l|c|c|c|c|c|c|}
        \hline
        \multirow{2}{*}{Benchmark Problem}
            & \multicolumn{2}{c|}{iLQR}
            & \multicolumn{2}{c|}{DDP}
            & \multicolumn{2}{c|}{Newton-LQR} \\
        \cline{2-7}
            & $\gamma$ & $N$ & $\gamma$ & $N$ & $\gamma$ & $N$ \\
        \hline
        Pendulum Swing-up & 0.656   & 0  & -11.86 & 2 & -4.87 & 1 \\
        \hline
        Cartpole Swing-up & 0.1272   & 0  & -1338.4 & 12 & -9.36 & 1 \\
        \hline
        Test Problem-1    & 0.02   & 0 & -37.07 & 4 & -36.69 & 1 \\
        \hline
        Test Problem-2    & 2$\times 10^{-4}$   & 0 &-0.0743 & 2486 & -0.0663 & 1  \\
        \hline
        Test Problem-3    & 2.88   & 0 & -41.619 & 49 & 2.88 & 0 \\
        \hline
        Test Problem-4    & -5.821$\times 10^{-14}$   & 49 & -404.17 & 97 & -0.28155 & 50 \\
        \hline
    \end{tabular}
\end{table}


Table~\ref{tab:Quu_PD} highlights the need for regularization in algorithms such as DDP and Newton-LQR, which incorporate second-order derivatives of the dynamics. Here, $\gamma$ denotes the minimum eigenvalue of $Q_{\vu\vu}^{t}$ over the entire time horizon, while $N$ represents the number of time steps at which $Q_{\vu\vu}^{t}$ fails to be positive definite. All results are evaluated at the first iteration of optimization.

For the pendulum swing-up, cart-pole swing-up, TP-1, and TP-2, the optimization is initialized with randomly generated control trajectories. For TP-3, the initial control sequence is chosen to alternate in sign between consecutive time steps, with $\vu_t=0.01\mathbf{1}_m$ for odd $t$ and $\vu_t=-0.01\mathbf{1}_m$ for even $t$. The problem parameters are $n=100$, $m=50$, $T=100$, and $\mu=1/75$. For TP-4, the initial control sequence corresponds to Starting Point 4 in Table~\ref{tab:starting_points}, with $n=100$, $m=10$, and $T=100$.

The results reveal that DDP frequently encounters a loss of positive definiteness, as indicated by negative values of $\gamma$ and a large number of non-positive-definite time steps. Newton-LQR exhibits similar behavior, further emphasizing the need for regularization when second-order dynamics information is incorporated. In contrast, iLQR maintains positive definiteness of $Q_{\vu\vu}^{t}$ throughout the horizon for all benchmark problems except TP-4. This observation is consistent with Proposition~\ref{prop:ilqr_descent}, as all problems except TP-4 satisfy the sufficient conditions for positive definiteness of $Q_{\vu\vu}^{t}$ and a predicted descent direction in iLQR.

It is important to note that the non-positive definiteness of certain $Q_{\vu\vu}^{t}$ matrices does not necessarily imply that the algorithm is unable to predict a cost descent direction. Table~\ref{tab:regularization_Req} summarizes, for 100 random initializations on the first four benchmark problems, the number of times each algorithm required regularization to ensure a descent direction and the maximum regularization value applied in these instances, all evaluated during the first iteration of the optimization process. As observed, iLQR does not necessitate any regularization, as it meets the assumptions stated in Proposition~\ref{prop:ilqr_descent}. In contrast, Newton-LQR requires some regularization, while DDP typically demands the largest regularization among the three algorithms.

\begin{table}[!htbp]
    \centering
    \caption{Frequency and magnitude of required regularization for different algorithms in ensuring a descent direction. The variable $N$ denotes the number of instances in which each algorithm necessitates regularization of $Q_{\vu\vu}^{t}$, while $\zeta_{\max}$ indicates the maximum regularization value applied across 100 random initializations of the control input during the first iteration of the optimization process.}
    \begin{tabular}{|l|c|c|c|c|c|c|c|c|}
    \hline
        \multirow{2}{*}{Algorithm} & \multicolumn{2}{|c|}{Pendulum} & \multicolumn{2}{c|}{Cart-pole} & \multicolumn{2}{c|}{TP-1} & \multicolumn{2}{c|}{TP-2} \\
        \cline{2-9} & $N$ & $\zeta_{\text{max}}$ &  $N$ & $\zeta_{\text{max}}$ &  $N$ & $\zeta_{\text{max}}$ &  $N$ & $\zeta_{\text{max}}$ \\
         \hline
         iLQR & 0 & 0 & 0 & 0 & 0 & 0 & 0 & 0\\
         \hline
         DDP & 6 & 2000 & 28 & 2000 & 1 & 0.02 & 0 & 0\\
         \hline
         Newton-LQR & 6 & 50 & 16 & 200 & 0 & 0 & 0 & 0\\
         \hline
    \end{tabular}
    
    \label{tab:regularization_Req}
\end{table}

To compare the performance of the algorithms under regularization, we evaluated all three methods across all six benchmark problems using the adaptive regularization strategy described in Section~\ref{sec:adap_reg}. For all experiments, the initial guess for the control inputs was set to zero. The convergence results are summarized in Table~\ref{tab:iterations_to_optimal}. For test problem 3, we used $n=100$, $m=50$, $T=100$, and $\mu = 1/75$, while for test problem 4, the parameters were $n=100$, $m=10$, and $T=100$. Across the six benchmarks, iLQR requires the fewest iterations or ties for the fewest iterations in every case. The relative performance of DDP and Newton-LQR is problem-dependent: Newton-LQR converges faster on the pendulum and TP-3/TP-4 benchmarks, whereas DDP requires fewer iterations on the cart-pole and TP-1/TP-2 benchmarks. For cart-pole swing-up task, DDP converged to a slightly higher cost than the optimal showing heavy regularization slowing down convergence. These results suggest that, even when regularization schemes are employed, the inclusion of Hessian terms does not necessarily improve convergence properties and, in some cases, may be inadvisable.




\begin{table}[htbp]
    \centering
    \caption{Number of iterations required for convergence to the optimal solution across all benchmark problems. The initial control sequence for each case was initialized to zero. $\cJ^{*}$ denotes the corresponding optimal cost.}
    \label{tab:test_problem_3_results_alt}
    \begin{tabular}{|l|c|c|c|c|c|c|}
        \hline
        \multirow{2}{*}{Problem}& \multicolumn{2}{|c|}{iLQR} &  \multicolumn{2}{c|}{DDP} &  \multicolumn{2}{c|}{Newton-LQR} \\
        \cline{2-7} & iter. & $\cJ_{\text{final}}$ & iter. & $\cJ_{\text{final}}$& iter. & $\cJ_{\text{final}}$\\
        \hline
        Pendulum         & 19 & 23.66       & 42  & 23.66     & 30 & 23.66\\
         Cart-pole         & 174 & 188.64        & 260 & \textbf{193.63}       & 323 & 188.64\\
          TP-1         & 8 & 6.7189        & 14 & 6.7189       & 18 & 6.7189 \\
           TP-2         & 16 & 41.6        & 19       & 41.6 & 27 & 41.6 \\
            TP-3         & 7 & 300.67       & 8  & 300.67      & 7 & 300.67 \\
             TP-4         & 4 & 8.5175        & 5 & 8.5175       & 4 & 8.5175 \\
        \hline
        
    \end{tabular}
\end{table}

Up to this point, our results indicate that incorporating the dynamics Hessian, as done in DDP and Newton-LQR, is generally not advantageous for the benchmark problems considered, since iLQR consistently attains faster convergence.  Our key observation is that the second-order approximation employed by DDP distorts the descent direction when the iterates are far from the optimum. A similar pattern is observed for Newton-LQR, which also exhibits slower convergence away from the optimal solution. In contrast, iLQR provides a guarantee that, by selecting a sufficiently small line search parameter, $\alpha$, the effect of nonlinearity in the dynamics can be minimized, thereby yielding more consistent descent directions.


In the next section, we will demonstrate that in the vicinity of the optimal solution, activating the dynamics Hessian can be beneficial to avoid prolonged convergence tails. In particular, we show that Newton-LQR surpasses DDP near the optimum, as it is a closer realization of the true Newton method relative to DDP, which can be regarded as a ``modified'' Newton approach.

\subsection{Faster Local Convergence Using Dynamics Hessians}
\label{subsec:hybrid_newton_results}

Equation~\eqref{eq:con_map} recalls the corresponding characterization of the convergence-rate. We now investigate whether these local convergence properties translate into practical computational advantages. However, as demonstrated in the preceding sections, iLQR can exhibit more favorable convergence behavior when initialized far from an optimum. Conversely, the additional curvature information incorporated by DDP and Newton-LQR can accelerate convergence in the neighborhood of a stationary solution. These observations motivate a hybrid optimization like that discussed in Section~\ref{subsec:hybrid_newton}.

To investigate the effectiveness of such a strategy, we consider the pendulum and cart-pole swing-up problems. We first obtain a sequence of nominal trajectories by running iLQR until convergence. At each iteration $k$, the corresponding iLQR trajectory is used to independently initialize DDP and Newton-LQR. We then record the number of additional iterations required by each algorithm to satisfy the same convergence criterion. This procedure allows us to identify the stage of optimization at which incorporating the dynamics Hessian information becomes advantageous.

\begin{figure*}[!t]
    \centering

    \begin{subfigure}[b]{0.48\textwidth}
        \centering
        \includegraphics[width=\linewidth]{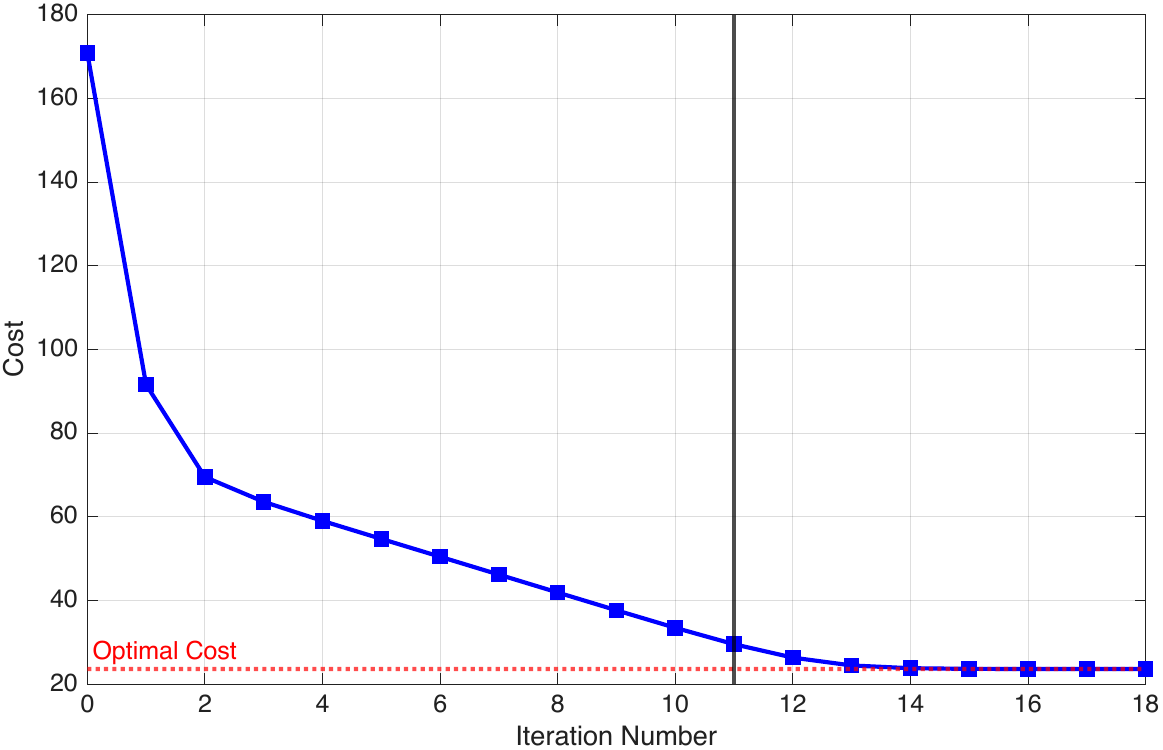}
        \caption{iLQR cost convergence.}
        \label{fig:ilqr_cost_convergence}
    \end{subfigure}
    \hfill
    \begin{subfigure}[b]{0.48\textwidth}
        \centering
        \includegraphics[width=\linewidth]{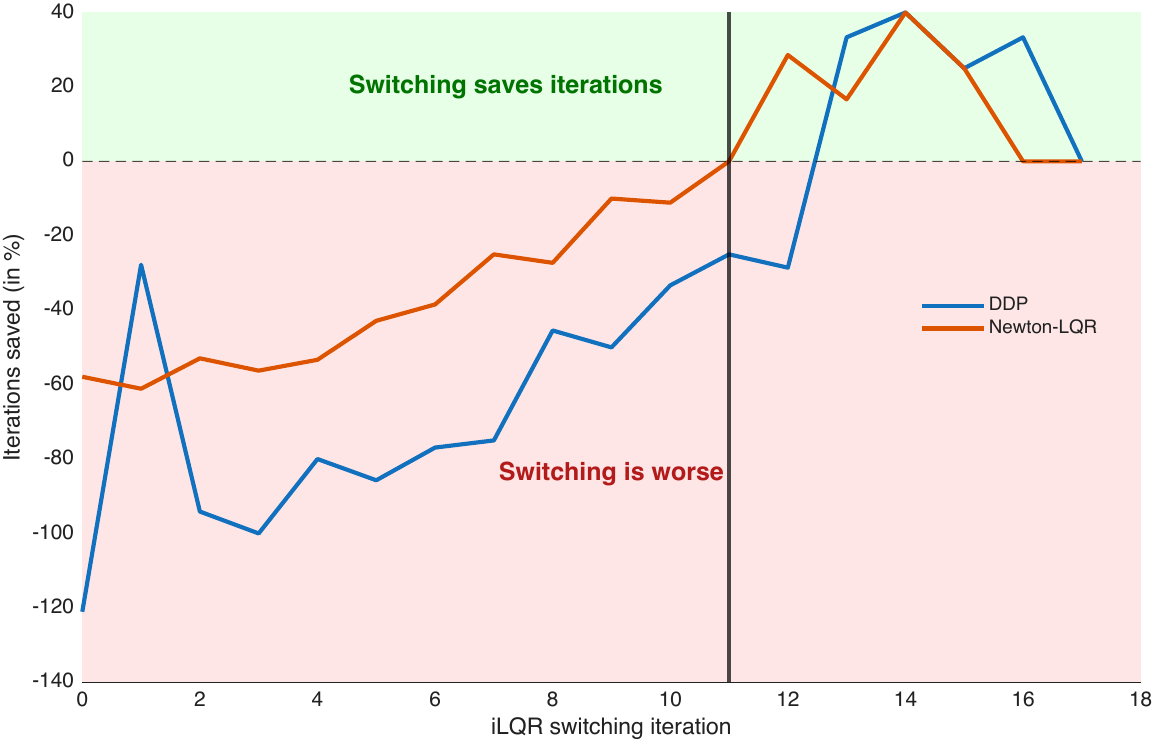}
        \caption{Percentage of iterations saved by switching.}
        \label{fig:remaining_iterations}
    \end{subfigure}

    \caption{Convergence comparison for the pendulum swing-up problem.
(a) Evolution of the objective function under iLQR, with the optimal
cost indicated by the red dotted line.
(b) Percentage of iLQR iterations saved by switching to DDP or
Newton-LQR at iteration $k$. The percentage is computed relative
to the remaining iLQR iterations.
Positive values (green region) indicate a reduction in iterations,
whereas negative values (red region) indicate that switching
requires more iterations than continuing with iLQR.
The vertical black line indicates the switching iteration $k=11$.}
    \label{fig:hybrid_convergence}
\end{figure*}

\begin{figure*}[!t]
    \centering

    \begin{subfigure}[b]{0.48\textwidth}
        \centering
        \includegraphics[width=\linewidth]{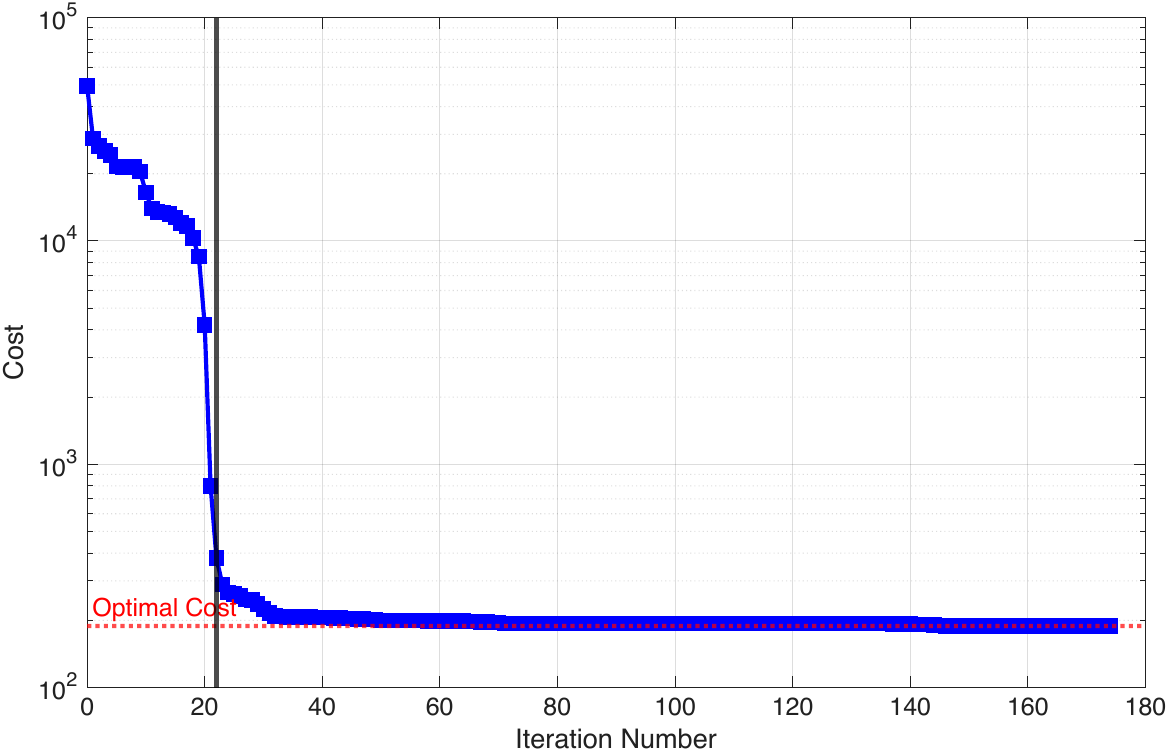}
        \caption{iLQR cost convergence.}
        \label{fig:ilqr_cost_cartpole}
    \end{subfigure}
    \hfill
    \begin{subfigure}[b]{0.48\textwidth}
        \centering
        \includegraphics[width=\linewidth]{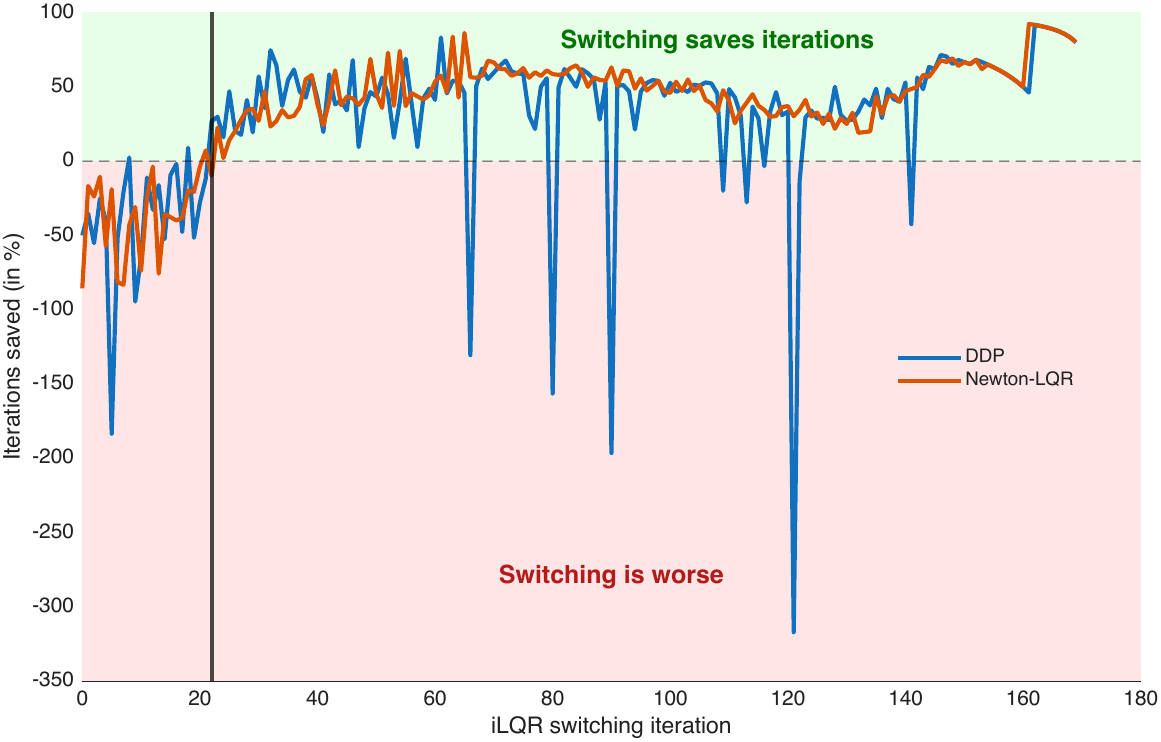}
        \caption{Percentage of iterations saved by switching.}
        \label{fig:remaining_iterations_cartpole}
    \end{subfigure}

    \caption{Convergence comparison for the cart-pole swing-up problem.
(a) Evolution of the objective function under iLQR, with the optimal
cost indicated by the red dotted line.
(b) Percentage of iLQR iterations saved by switching to DDP or
Newton-LQR at iteration $k$. The percentage is computed relative
to the remaining iLQR iterations.
Positive values (green region) indicate a reduction in iterations,
whereas negative values (red region) indicate that switching
requires more iterations than continuing with iLQR.
The vertical black line indicates the switching iteration $k=22$.}
    \label{fig:hybrid_convergence_cartpole}
\end{figure*}

Figure~\ref{fig:hybrid_convergence} presents the convergence results for the pendulum swing-up problem. Figure~\ref{fig:ilqr_cost_convergence} illustrates the reduction in the objective function over successive iLQR iterations, while Figure~\ref{fig:remaining_iterations} compares the additional iterations required by DDP and Newton-LQR to converge relative to iLQR. Specifically, the relative difference in the number of iterations is defined as
\begin{align}
    \Delta_{\mathrm{iter}}(\%)
    = \frac{N_{\mathrm{iLQR}}- N_{\mathrm{algo}}}
    {N_{\mathrm{iLQR}}} \times 100,
\end{align}
where $N_{\mathrm{algo}}$ denotes the number of remaining iterations required by DDP or Newton-LQR, and $N_{\mathrm{iLQR}}$ denotes the corresponding number for iLQR. The zero baseline, represented by a dashed line, indicates identical convergence performance. The green region ($\Delta_{\mathrm{iter}}>0$) identifies iterations at which switching to a second-order method reduces the number of iterations required for convergence, whereas the red region ($\Delta_{\mathrm{iter}}<0$) indicates that switching is less effective than continuing with iLQR.

For the pendulum swing-up problem, Newton-LQR generally requires fewer iterations than DDP, although their performance becomes nearly identical toward the end of optimization. Notably, Newton-LQR crosses below the iLQR baseline at approximately the 11th iteration, followed by DDP two iterations later. A solid black vertical line marks the 11th iteration in both figures, highlighting the onset of improved convergence when dynamics-Hessian information is incorporated. This observation is consistent with Figure~\ref{fig:ilqr_cost_convergence}, where the objective function has already approached its final value by approximately the 11th iteration. Thus, while iLQR exhibits favorable convergence during the initial stages of optimization, incorporating second-order dynamics information becomes beneficial as the nominal trajectory approaches the optimum.

A similar trend is observed for the cart-pole swing-up problem, as shown in Figure~\ref{fig:hybrid_convergence_cartpole}. Figure~\ref{fig:ilqr_cost_cartpole} illustrates the progressive reduction in the objective function under iLQR, while Figure~\ref{fig:remaining_iterations_cartpole} shows that both DDP and Newton-LQR begin to outperform iLQR after approximately 22 iterations. This transition is marked by a solid black vertical line in the corresponding figures and coincides with the stage at which the objective function has approached its final value.

However, an important distinction emerges between the two second-order methods. While Newton-LQR consistently remains in the green region beyond the switching point, indicating a reduction in the number of remaining iterations, DDP exhibits several spikes that occasionally exceed the iLQR baseline. These fluctuations suggest that DDP is more sensitive to the nominal trajectory, whereas Newton-LQR provides more consistent improvements in convergence. This behavior is consistent with the analysis in Section~\ref{sec:ipfocs}, which establishes that DDP approximates the Newton-LQR formulation, with the approximation becoming increasingly accurate as the trajectory approaches stationarity.

The numerical results highlight two important characteristics of these algorithms. First, incorporating the dynamics Hessian information does not necessarily improve convergence when the iterates are far from an optimum, where iLQR can provide more reliable progress. Second, as the iterates approach a stationary solution, the additional second-order information becomes increasingly beneficial, allowing DDP and Newton-LQR to converge in fewer iterations. Among the two methods, Newton-LQR exhibits more consistent local convergence behavior across the examples considered. These findings motivate the proposed hybrid approach, in which iLQR is employed initially to approach the optimal solution, followed by Newton-LQR to exploit its faster local convergence. Such a strategy seeks to combine the favorable early-stage convergence behavior of iLQR with the local efficiency of Newton's method. 

\section{Conclusion}\label{sec:conclusion}
In this paper, we have investigated discrete-time nonlinear optimal control problems by analyzing Newton's method, Differential Dynamic Programming (DDP), and the iterative Linear Quadratic Regulator (iLQR) through a unified Sequential Quadratic Programming (SQP) framework. The central insight is that these algorithms solve a sequence of Quadratic Programs (QPs), which reduce to Linear Quadratic Regulator (LQR) problems by exploiting the structure of optimal control. Within this framework, we introduce ``Newton-LQR'', a novel Riccati-based formulation for solving the exact Newton-SQP quadratic subproblem, and establish its relationship with iLQR and DDP.

Our analysis reveals that DDP can be interpreted as an approximation of Newton-LQR, with their formulations agreeing near a stationary solution under appropriate conditions. However, the inclusion of dynamics Hessians can render the local quadratic subproblems indefinite, necessitating regularization and potentially producing unreliable descent directions. In contrast, iLQR provides consistent descent directions under suitable convexity assumptions and admits theoretical convergence guarantees.

Numerical experiments across six benchmark problems corroborate the theoretical findings, demonstrating that iLQR generally requires less regularization and converges faster than DDP and Newton-LQR. We further establish that heavy regularization causes the control update to asymptotically approach a scaled gradient-descent step, potentially slowing convergence. Nevertheless, our numerical investigations demonstrate that incorporating dynamics Hessians can accelerate convergence near a stationary solution, where Newton-LQR exhibits more consistent improvements than DDP in the swing-up examples considered.

These findings motivate a hybrid strategy that employs iLQR during the initial optimization phase and transitions to Newton-LQR near optimality. Overall, our results establish iLQR as a principled SQP approach rather than merely an approximation of DDP, while highlighting the importance of balancing reliable descent directions with the local convergence advantages of second-order information.


\bibliographystyle{plain}        
\bibliography{references}           



\newpage
\appendix
\begin{center}
    \Large APPENDIX
\end{center}

\section{Impact of Regularization on Predicted Cost Reduction}
\label{app:regularization}

Equation~\eqref{eq:exp_cost_red} gives the predicted change in
cost in the absence of regularization. When regularization is
introduced, the predicted cost reduction takes the form
\begin{align}
    \Delta \cJ_{\mathrm{pred}}
    &=
    \alpha\sum_{t=0}^{T-1}
    (Q_{\vu}^{t})^{\top}\vk_t
    +
    \frac{\alpha^2}{2}
    \sum_{t=0}^{T-1}
    \vk_t^{\top}Q_{\vu\vu}^{t}\vk_t,
    \label{eq:app_exp_cost_reg}
\end{align}
where $\alpha\in(0,1]$ is the line-search parameter and
\begin{align}
    \vk_t
    =
    -(Q_{\vu\vu}^{t}+\zeta I)^{-1}Q_{\vu}^{t}.
\end{align}

Substituting the regularized feedforward gain gives
\begin{align}
    \Delta \cJ_{\mathrm{pred}}
    &=
    -\alpha\sum_{t=0}^{T-1}
    (Q_{\vu}^{t})^{\top}
    (Q_{\vu\vu}^{t}+\zeta I)^{-1}
    Q_{\vu}^{t}
    \nonumber\\
    &\quad+
    \frac{\alpha^2}{2}
    \sum_{t=0}^{T-1}
    (Q_{\vu}^{t})^{\top}
    (Q_{\vu\vu}^{t}+\zeta I)^{-1}
    Q_{\vu\vu}^{t}
    (Q_{\vu\vu}^{t}+\zeta I)^{-1}
    Q_{\vu}^{t}.
\end{align}

Owing to the symmetry of $Q_{\vu\vu}^{t}$, consider its
eigendecomposition
\begin{align}
    Q_{\vu\vu}^{t}
    =
    V_t\Lambda_t V_t^{\top},
\end{align}
where $V_t$ is an orthogonal matrix containing the
eigenvectors and
\begin{align}
    \Lambda_t
    =
    \operatorname{diag}
    (\lambda_1^t,\lambda_2^t,\ldots,\lambda_m^t)
\end{align}
is the diagonal matrix of corresponding eigenvalues.

Using this decomposition, the predicted cost reduction
can be expressed as
\begin{align}
    \Delta \cJ_{\mathrm{pred}}
    &=
    -\alpha\sum_{t=0}^{T-1}
    (Q_{\vu}^{t})^{\top}
    V_t(\Lambda_t+\zeta I)^{-1}
    V_t^{\top}Q_{\vu}^{t}
    \nonumber\\
    &\quad+
    \frac{\alpha^2}{2}
    \sum_{t=0}^{T-1}
    (Q_{\vu}^{t})^{\top}V_t
    (\Lambda_t+\zeta I)^{-1}
    \Lambda_t
    (\Lambda_t+\zeta I)^{-1}
    V_t^{\top}Q_{\vu}^{t}.
\end{align}

Collecting the terms yields
\begin{align}
    \Delta \cJ_{\mathrm{pred}}
    &=
    \sum_{t=0}^{T-1}
    (Q_{\vu}^{t})^{\top}V_t
    \operatorname{diag}
    \left[
    \frac{
    (\alpha^2-2\alpha)\lambda_i^t-2\alpha\zeta
    }{
    2(\lambda_i^t+\zeta)^2
    }
    \right]
    V_t^{\top}Q_{\vu}^{t}.
    \label{eq:app_eigen_cost_red}
\end{align}

The coefficient associated with the $i$-th eigenmode
at time $t$ is
\begin{align}
    \eta_i^t
    &=
    \frac{
    (\alpha^2-2\alpha)\lambda_i^t-2\alpha\zeta
    }{
    2(\lambda_i^t+\zeta)^2
    }.
\end{align}

For $\alpha\in(0,1]$ and $\zeta\geq0$, this coefficient
is strictly negative whenever
\begin{align}
    \zeta>\left( \frac{\alpha}{2} - 1\right) \lambda_{i}^{t}.
\end{align}

Therefore, if $Q_{\vu\vu}^{t}$ is already positive
definite, no regularization is required to ensure a
negative predicted cost change. However, when
$Q_{\vu\vu}^{t}$ possesses non-positive eigenvalues,
the regularization parameter can be chosen such that
\begin{align}
    \boxed{
    \zeta>
    -\min_{t}
    \lambda_{\min}(Q_{\vu\vu}^{t}),
    \qquad \zeta\geq0.
    }
    \label{eq:app_reg_condition}
\end{align}

This ensures that
\begin{align}
    \boxed{
    \Delta\cJ_{\mathrm{pred}}<0
    }
\end{align}
for every $\alpha\in(0,1]$, provided that the nominal is not a stationary point($Q_{\vu}^{t}\neq \bo$).

Thus, sufficiently large regularization guarantees
a descent direction for the local quadratic prediction.
A corresponding reduction in the actual nonlinear
cost must be verified through the forward-pass
line-search procedure.

\subsection{Heavy Regularization and Gradient Descent}

We now consider the limiting case of heavy
regularization, $\zeta\gg0$, and investigate its
effect on the control update.

The feedforward gain is given by
\begin{align}
    \vk_t
    &=
    -(Q_{\vu\vu}^{t}+\zeta I)^{-1}Q_{\vu}^{t}
    \nonumber\\
    &=
    -V_t(\Lambda_t+\zeta I)^{-1}
    V_t^{\top}Q_{\vu}^{t}.
\end{align}

For sufficiently large $\zeta$, we have
\begin{align}
    \frac{1}{\lambda_i^t+\zeta}
    =
    \frac{1}{\zeta}
    +
    \mathcal O\left(\frac{1}{\zeta^2}\right).
\end{align}

Consequently,
\begin{align}
    \vk_t
    &=
    -\frac{1}{\zeta}
    V_tV_t^{\top}Q_{\vu}^{t}
    +
    \mathcal O\left(\frac{1}{\zeta^2}\right).
\end{align}

Since $V_tV_t^{\top}=I$, we obtain
\begin{align}
    \boxed{
    \vk_t
    =
    -\frac{1}{\zeta}Q_{\vu}^{t}
    +
    \mathcal O\left(\frac{1}{\zeta^2}\right).
    }
    \label{eq:app_heavy_k}
\end{align}

Similarly, the feedback gain satisfies
\begin{align}
    \vK_t
    &=
    -(Q_{\vu\vu}^{t}+\zeta I)^{-1}
    Q_{\vu\vx}^{t}
    \nonumber\\
    &=
    -V_t(\Lambda_t+\zeta I)^{-1}
    V_t^{\top}Q_{\vu\vx}^{t}.
\end{align}

Hence,
\begin{align}
    \boxed{
    \vK_t
    =
    -\frac{1}{\zeta}Q_{\vu\vx}^{t}
    +
    \mathcal O\left(\frac{1}{\zeta^2}\right).
    }
    \label{eq:app_heavy_K}
\end{align}

The forward-pass control update is
\begin{align}
    \delta\vu_t
    =
    \alpha\vk_t+\vK_t\delta\vx_t.
\end{align}

From eq.~\eqref{eq:lin_exp_dyn}, state perturbation is a function $\vk_{t}$ and $\vK_{t}$. Since both $\vk_t$ and $\vK_t$ are
$\mathcal O(\zeta^{-1})$, it follows by induction,
for a fixed finite horizon and bounded dynamics
Jacobians, that
\begin{align}
    \delta\vx_t
    =
    \mathcal O\left(\frac{1}{\zeta}\right).
\end{align}

Therefore, the feedback correction satisfies
\begin{align}
    \boxed{
    \vK_t\delta\vx_t
    =
    \mathcal O\left(\frac{1}{\zeta^2}\right).
    }
\end{align}

Substituting these results into the control update,
we obtain
\begin{align}
    \delta\vu_t
    &=
    \alpha\vk_t
    +
    \mathcal O\left(\frac{1}{\zeta^2}\right)
    \nonumber\\
    &=
    -\frac{\alpha}{\zeta}Q_{\vu}^{t}
    +
    \mathcal O\left(\frac{1}{\zeta^2}\right).
    \label{eq:app_heavy_control}
\end{align}

Furthermore,
\begin{align}
    Q_{\vu}^{t}
    =
    \frac{\partial\cJ}{\partial\vu_t},
\end{align}
where the derivative is taken with respect to the
control sequence after eliminating the states
through the nonlinear dynamics. 

Thus,
\begin{align}
    \boxed{
    \delta\vu_t
    =
    -\frac{\alpha}{\zeta}
    \frac{\partial\cJ}{\partial\vu_t}
    +
    \mathcal O\left(\frac{1}{\zeta^2}\right).
    }
    \label{eq:app_gradient_limit}
\end{align}
This is precisely the leading-order form of a
gradient-descent update, with an effective step
size $\alpha/\zeta$.

\begin{proposition}
\label{prop:heavy_reg}
Under assumption~\ref{assum:smoothness}, under heavy
regularization,  the control updates generated
by iterative perturbation feedback methods asymptotically
reduce to gradient descent:
\begin{align}
    \delta\vu_t
    =
    -\frac{\alpha}{\zeta}
    \frac{\partial\cJ}{\partial\vu_t}
    +
    \mathcal O\left(\frac{1}{\zeta^2}\right).
\end{align}
\end{proposition}

This result demonstrates that excessive
regularization progressively suppresses the
second-order information incorporated by DDP
and Newton-LQR. Consequently, although
regularization can ensure a well-posed
quadratic subproblem and a predicted descent
direction, large regularization
reduces the leading-order control update to
a scaled gradient step and can substantially
slow convergence.

\section{Deriving DDP, Stagewise-Newton, and iLQR using Riccati equation}
 In the following, we reiterate how Differential Dynamic Programming (DDP), Stagewise-Newton method, and the iterative Linear Quadratic Regulator (iLQR) utilize the sparsity structure inherent in the optimal control problem described by eq.~\eqref{eq:ocp} and eq.~\eqref{eq:dynamics}. This exploitation of sparsity leads to significantly reduced computational complexity compared to conventional Sequential Quadratic Programming (SQP).

To begin with, consider a feasible control sequence and a corresponding trajectory $(\bar{\vx}_t, \bar{\vu}_t)$ can be constructed by evolving the dynamics.
Next, the augmented cost (using eqs. \eqref{eq:ocp} and \eqref{eq:dynamics}) can be written as:
\begin{equation}
    \min_{\vx_t, \vu_{t}, \blambda_t} \sum_{t=0}^{T-1}\Big( c(\vx_{t}, \vu_{t}) + \blambda_{t+1}\tran(\vf(\vx_{t}, \vu_{t}) - \vx_{t+1}) \Big) + C_{T}(\vx_{T})
\end{equation}
where $\blambda_{t+1} \in \mathbb{R}^{n}$ are the Lagrange multipliers/costates. 

In addition, let $(\bar{\vx}_t, \bar{\vu}_t, \bar{\blambda}_t)$ be the current guess to the problem. This augmented cost can be expanded about $(\bar{\vx}_{t}, \bar{\vu}_{t}, \bar{\blambda}_{t})$ to get
\begin{dmath}
\bar{\cJ}+\delta\cJ+\delta^{2}\cJ
=\min_{\delta \vx_t,\;\delta \vu_t,\;\delta \blambda_t}\Bigg(
\sum_{t=0}^{T-1}\Big(
c(\bar{\vx}_{t},\bar{\vu}_{t})
+\bar{\blambda}_{t+1}\tran\big(\vf(\bar{\vx}_{t},\bar{\vu}_{t})-\bar{\vx}_{t+1}\big)
\Big)
+ C_{T}(\bar{\vx}_{T}) 
+\sum_{t=0}^{T-1}\Big(
c_{\vx_{t}}\delta \vx_{t}
+c_{\vu_{t}}\delta \vu_{t}
+\delta \blambda_{t+1}\tran\big(\bar{\vx}_{t+1}-\vf(\bar{\vx}_{t},\bar{\vu}_{t})\big)
+\bar{\blambda}_{t+1}\tran\big(
\vf_{\vx_{t}}\delta \vx_{t}
+\vf_{\vu_{t}}\delta \vu_{t}
-\delta \vx_{t+1}
\big)
+\tfrac{1}{2}\delta \vx_{t}\tran c_{\vx_{t} \vx_{t}} \delta \vx_{t}
+\tfrac{1}{2}\delta \vu_{t}\tran c_{\vu_{t} \vu_{t}} \delta \vu_{t}
+\delta \vx_{t}\tran c_{\vx_{t}\vu_{t}} \delta \vu_{t}
+\bar{\blambda}_{t+1}\tran\Big(
\tfrac{1}{2}\delta \vx_{t}\tran \vf_{\vx_{t} \vx_{t}} \delta \vx_{t}
+\delta \vx_{t}\tran \vf_{\vx_{t} \vu_{t}} \delta \vu_{t}
+\tfrac{1}{2}\delta \vu_{t}\tran \vf_{\vu_{t} \vu_{t}} \delta \vu_{t}
\Big)
+\delta \blambda_{t+1}\tran\big(
\vf_{\vx_{t}}\delta \vx_{t}
+\vf_{\vu_{t}}\delta \vu_{t}
-\delta \vx_{t+1}
\big)
\Big)
+ C_{T,\vx}\delta \vx_{T}
+\tfrac{1}{2}\delta \vx_{T}' C_{T,\vx\vx} \delta \vx_{T}
\Bigg).
\label{eq:per_optim}
\end{dmath}
The necessary conditions for the above QP are:
\begin{subequations}
    \begin{dmath}
        \vf(\bar{\vx}_t,\bar{\vu}_t) + \vf_{\vx_t}\delta \vx_t + \vf_{\vu_t}\delta \vu_t - \delta \vx_{t + 1} - \bar{\vx}_{t+1} = \bo
        \label{eq:lin_dyn}
    \end{dmath}
    \begin{dmath}
        c_{\vx_t} + c_{\vx_t \vx_t}\delta \vx_t + c_{\vx_t \vu_t}\delta \vu_t
        +(\vf_{\vx_t \vx_t} \delta \vx_{t} + \vf_{\vx_t \vu_{t}}\delta \vu_{t})\tran\bar{\blambda}_{t+1} + \vf_{\vx_t}\tran(\bar{\blambda}_{t+1} + \delta \blambda_{t+1})
        -(\bar{\blambda}_t + \delta \blambda_t) = \bo
        \label{eq:state_DDP}
    \end{dmath}
    \begin{dmath}
        c_{\vu_t} + c_{\vu_t \vu_t}\delta \vu_t + c_{\vu_t \vx_t} \delta \vx_t
        +(\vf_{\vu_t \vu_t} \delta \vu_{t} + \vf_{\vu_t, \vx_{t}}\delta \vx_{t})\tran \bar{\blambda}_{t+1} + \vf_{\vu_t}\tran(\bar{\blambda}_{t+1} + \delta \blambda_{t+1}) = \bo
        \label{eq:con_DDP}
    \end{dmath}    
\end{subequations}
with the boundary/transversality condition:
\begin{align}
    \bar{\blambda}_{T} + \delta \blambda_{T} = C_{T,\vx} + C_{T, \vx\vx}\delta \vx_{T}, \nonumber \\
    \bar{\blambda}_{T} = C_{T,\vx} \quad \text{and} \quad \delta \blambda_{T} = C_{T,\vx\vx} \delta \vx_{T}.
    \label{eq:transver}
\end{align}
During the forward pass, the dynamics in Eq.~\eqref{eq:dynamics} are propagated while updating the control trajectory. In the backward pass, the associated multipliers and a linear feedback control law are computed implicitly. A suitable line-search parameter $\alpha$ is then selected to ensure feasibility of the updated trajectory. We now describe these algorithms in the context of the augmented perturbed-cost optimization problem formulated in Eq.~\eqref{eq:per_optim}.

\subsection{DDP}
We give an overview of these algorithms starting with ``DDP". First, it starts with an initial guess of $\{\vu_{t}\}_{t=0}^{T-1}$ and initializes the trajecotry using the dynamics $\bar{\vx}_{t+1} = \vf(\bar{\vx}_t,\bar{\vu}_t)$ in the forward pass. Following this, the control update $\delta \vu_t$ is solved by leveraging the sparsity of variables in the backward pass. The backward pass begins at the last time step  with the transversality condition in eq. \eqref{eq:transver} and then eq. \eqref{eq:con_DDP} is used to find the control update. Furthermore, $\delta \vx_{T}$ is replaced by linearized dynamics using eq. \eqref{eq:lin_dyn} in terms of $\delta \vu_{t-1}$  and $\delta \vx_{t-1}$. The equation which follows from \eqref{eq:con_DDP} at time-step $T-1$ is
\begin{dmath*}
    c_{\vu_{T-1}} + c_{\vu_{T-1}\vu_{T-1}}\delta \vu_{T-1} + c_{\vu_{T-1}\vx_{T-1}}\delta \vx_{T-1} + (\vf_{\vu_{T-1}\vu_{T-1}}\delta \vu_{T-1} + \vf_{\vu_{T-1}\vx_{T-1}}\delta \vx_{T-1})\tran\bar{\blambda}_{T} +
    \vf_{\vu_{T-1}}\tran(\bar{\blambda}_{T} + \delta \blambda_{T}) = \bo.
\end{dmath*}
Using the transversality condition in eq. \eqref{eq:transver}, 
\begin{dmath*}
    c_{\vu_{T-1}} + c_{\vu_{T-1}\vu_{T-1}}\delta \vu_{T-1} + c_{\vu_{T-1}\vx_{T-1}}\delta \vx_{T-1}+ (\vf_{\vu_{T-1}\vu_{T-1}}\delta \vu_{T-1} + \vf_{\vu_{T-1}\vx_{T-1}}\delta \vx_{T-1})\tran C_{T,\vx}
    +\vf_{\vu_{T-1}}\tran(C_{T,\vx} + C_{T,\vx\vx}\delta \vx_{T}) = \bo.
\end{dmath*}
Furthermore, using the linearized dynamics in eq. \eqref{eq:lin_dyn} and replacing $\delta \vx_{T}$ in the above equation, we get
\begin{dmath*}
    c_{\vu_{T-1}} + c_{\vu_{T-1}\vu_{T-1}}\delta \vu_{T-1} + c_{\vu_{T-1}\vx_{T-1}}\delta \vx_{T-1} 
    + (\vf_{\vu_{T-1}\vu_{T-1}}\delta \vu_{T-1} + \vf_{\vu_{T-1}\vx_{T-1}}\delta \vx_{T-1})\tran C_{T,\vx}
    +\vf_{\vu_{T-1}}'(C_{T,\vx} + C_{T,\vx\vx}(\vf_{\vx_{T-1}} \delta \vx_{T-1} + \vf_{\vu_{T-1}}\delta \vu_{T-1})) = \bo.
\end{dmath*}
The above equation is linear in terms of $\delta \vu_{T-1}$ and the solution to it is:
\begin{dmath}
    \delta \vu_{T-1} = 
    -\Big{(} c_{\vu_{T-1}\vu_{T-1}} + \vf_{\vu_{T-1}\vu_{T-1}} \otimes C_{T,\vx} + \vf_{\vu_{T-1}}\tran C_{T,\vx\vx}\vf_{\vu_{T-1}}\Big{)}^{-1}
    \Bigg{(}c_{\vu_{T-1}} + \vf_{\vu_{T-1}}\tran C_{T,x} +
    \Big{(} c_{\vu_{T-1}\vx_{T-1}} + \vf_{\vu_{T-1}\vx_{T-1}} \otimes C_{T,\vx} + \vf_{\vu_{T-1}}C_{T,\vx\vx}\vf_{\vx_{T-1}}\Big{)}\delta \vx_{T-1}\Bigg{)}.
\end{dmath}
The control update in the above equation can be written in short as
\begin{align}
    \delta \vu_{T-1} = -(Q_{\vu\vu}^{T-1})^{-1}(Q_{\vu}^{T-1} + Q_{\vu\vx}^{T-1}\delta \vx_{T-1}).
    \label{eq:con_terminal}
\end{align}
The above term is the linear feedback for the DDP. The terms $Q_{\vu\vu}$,$Q_{\vu\vx}$ and $Q_{\vu}$ are
\begin{align*}
    &Q_{\vu\vu}^{T-1} = c_{\vu_{T-1}\vu_{T-1}} + \vf_{\vu_{T-1}\vu_{T-1}} \otimes C_{T,\vx} + \vf_{\vu_{T-1}}\tran C_{T,\vx\vx}\vf_{\vu_{T-1}}.\\
    &Q_{\vu\vx}^{T-1} = c_{\vu_{T-1}\vx_{T-1}} + \vf_{\vu_{T-1}\vx_{T-1}} \otimes C_{T,\vx} + \vf_{\vu_{T-1}}C_{T,\vx\vx}\vf_{\vx_{T-1}}.\\
   &Q_{\vu}^{T-1} = c_{\vu_{T-1}} + \vf_{\vu_{T-1}}\tran C_{T,\vx}.
\end{align*}

Hence, eq. \eqref{eq:con_DDP} is solved to get a linearized feedback control update at the terminal step, $\delta \vu_{T-1} = \vk_{T-1} + \vK_{T-1}\delta \vx_{T-1}$. This update is used in eq. \eqref{eq:state_DDP}, along with linearized dynamics from eq. \eqref{eq:lin_dyn}, to get a backward pass on $\bar{\blambda}_{T}$ (represented as $\bv_T$) and $\mathbf{S}_{T}$, where $\delta \blambda_{T} = \mathbf{S}_{T} \delta \vx_{T}$.
It comes from eq. \eqref{eq:state_DDP}:
\begin{dmath*}
    c_{\vx_{T-1}} + c_{\vx_{T-1}\vx_{T-1}}\delta \vx_{T-1} + c_{\vx_{T-1}\vu_{T-1}}\delta \vu_{T-1}
    +(\vf_{\vx_{T-1}\vx_{T-1}}\delta \vx_{T-1}+ \vf_{\vx_{T-1}\vu_{T-1}}\delta \vu_{T-1})\tran \bar{\blambda}_{T} + \vf_{\vx_{T-1}}\tran (\bar{\blambda}_{T} + \delta \blambda_{T})
    -(\bar{\blambda}_{T-1} + \delta \blambda_{T-1}) = \boldsymbol{0}.
\end{dmath*}
Substituting the value of control $\delta \vu_{T-1}$ from eq. \eqref{eq:con_terminal} and using the transversality condition one gets
\begin{dmath*}
    c_{\vx_{T-1}} + c_{\vx_{T-1}\vx_{T-1}}\delta \vx_{T-1} + c_{\vx_{T-1}\vu_{T-1}}(-Q_{\vu\vu}^{-1}(Q_{\vu}^{T-1} + Q_{\vu\vx}^{T-1}\delta \vx_{T-1})) 
    +(\vf_{\vx_{T-1}\vx_{T-1}}\delta \vx_{T-1}+ \vf_{\vx_{T-1}\vu_{T-1}}(-Q_{\vu\vu}^{-1}(Q_{\vu}^{T-1} + Q_{\vu\vx}^{T-1}\delta \vx_{T-1})))\tran C_{T,\vx} 
    +\vf_{\vx_{T-1}}\tran (C_{T,\vx} + C_{T,\vx\vx}\delta \vx_{T})
    -(\bar{\blambda}_{T-1} + \delta \blambda_{T-1}) = \bo.
\end{dmath*}
Using linearized dynamics for $\delta \vx_{T}$ from eq. \eqref{eq:lin_dyn} and replacing the control from eq. \eqref{eq:con_terminal}, the above expression becomes
\begin{dmath*}
    c_{\vx_{T-1}} + c_{\vx_{T-1}\vx_{T-1}}\delta \vx_{T-1} + c_{\vx_{T-1}\vu_{T-1}}(-(Q_{\vu\vu}^{T-1})^{-1}(Q_{\vu}^{T-1} + Q_{\vu\vx}^{T-1}\delta \vx_{T-1})) 
    +(\vf_{\vx_{T-1},x_{T-1}}\delta \vx_{T-1}+ \vf_{\vx_{T-1}\vu_{T-1}}(-(Q_{\vu\vu}^{T-1})^{-1}(Q_{\vu}^{T-1} + Q_{\vu\vx}^{T-1}\delta \vx_{T-1})))\tran C_{T,\vx} 
    + \vf_{\vx_{T-1}}\tran (C_{T,\vx} + C_{T,\vx\vx}(\vf_{\vx_{T-1}}\delta \vx_{T-1} + \vf_{\vu_{T-1}}
    (-(Q_{\vu\vu}^{T-1})^{-1}(Q_{\vu}^{T-1} + Q_{\vu\vx}^{T-1}\delta \vx_{T-1}))))
    -(\bar{\blambda}_{T-1} + \delta \blambda_{T-1}) = \bo.
    \label{eq:back_pass_DDP1}
\end{dmath*}
It is evident from the above equation that $\bar{\blambda}_{T-1}$ and $\delta \blambda_{T-1}$ can be separated using the form $\bv_{T-1} + \vS_{T-1}\delta \vx_{T-1}$. Separating out the constant terms gives
\begin{dmath*}
    c_{\vx_{T-1}} - c_{\vx_{T-1}\vu_{T-1}}(Q_{\vu\vu}^{T-1})^{-1}Q_{\vu}^{T-1} - (\vf_{\vx_{T-1}\vu_{T-1}}(Q_{\vu\vu}^{T-1})^{-1}Q_{\vu}^{T-1})\tran C_{T,\vx}
    +\vf_{\vx_{T-1}}\tran C_{T,\vx} - \vf_{\vx_{T-1}}\tran C_{T,\vx\vx}\vf_{\vu_{T-1}}(Q_{\vu\vu}^{T-1})^{-1}Q_{\vu}^{T-1} - \bv_{T-1} = \bo.
\end{dmath*}
Rearranging the terms give,
\begin{dmath*}
    c_{\vx_{T-1}} + \vf_{\vx_{T-1}}\tran C_{T,\vx}
    -(c_{\vx_{T-1}\vu_{T-1}} + C_{T,\vx}\otimes \vf_{\vx_{T-1}\vu_{T-1}} + \vf_{\vx_{T-1}}\tran C_{T,\vx\vx}\vf_{\vu_{T-1}})(Q_{\vu\vu}^{T-1})^{-1}Q_{\vu}^{T-1} = \bv_{T-1},
\end{dmath*}
or in short,
\begin{dmath}
    \bv_{T-1} = Q_{\vx}^{T-1} - Q_{\vx\vu}^{T-1}(Q_{\vu\vu}^{T-1})^{-1}Q_{\vu}^{T-1},
\end{dmath}
where
\begin{align*}
    &Q_{\vx}^{T-1} = c_{\vx_{T-1}} + \vf_{\vx_{T-1}}\tran C_{T,\vx},\quad \text{and}\\
    &Q_{\vx\vu}^{T-1} =  c_{\vx_{T-1}\vu_{T-1}} +  C_{T,\vx}\otimes \vf_{\vx_{T-1}\vu_{T-1}}+ \vf_{\vx_{T-1}}\tran C_{T,\vx\vx}\vf_{\vu_{T-1}}.
\end{align*}
Similarly, the terms involving $\delta \vx_{T-1}$ in eq.~\eqref{eq:back_pass_DDP1} are:
\begin{dmath*}
    c_{\vx_{T-1}\vx_{T-1}} - c_{\vx_{T-1}\vu_{T-1}}(Q_{\vu\vu}^{T-1})^{-1}Q_{\vu\vx} + (\vf_{\vx_{T-1}\vx_{T-1}}\otimes C_{T,\vx})
    - (\vf_{\vx_{T-1}\vu_{T-1}}\otimes C_{T,\vx}  Q_{\vu\vu}^{-1}Q_{\vu\vx}) + \vf_{\vx_{T-1}}\tran C_{T,\vx\vx}\vf_{\vx_{T-1}} 
    - \vf_{\vx_{T-1}}\tran C_{T,\vx\vx}\vf_{\vu_{T-1}}Q_{\vu\vu}^{-1}Q_{\vu\vx} = \vS_{T-1},
\end{dmath*}
or in short,
\begin{align}
    \vS_{T-1} = Q_{\vx\vx}^{T-1} - Q_{\vx\vu}^{T-1}(Q_{\vu\vu}^{T-1})^{-1}Q_{\vu\vx}^{T-1}
\end{align}
where
\begin{dmath*}
    Q_{\vx\vx}^{T-1} = c_{\vx_{T-1}\vx_{T-1}} +  \vf_{\vx_{T-1} \vx_{T-1}} \otimes C_{T,\vx} +  \vf_{\vx_{T-1}}\tran C_{T,\vx\vx}\vf_{\vx_{T-1}}
\end{dmath*}
Now, this can be generalized for every step to get $\bv_{t}$ and $\vS_{t}$ and these can be used to calculate $\vu_{t} = \vk_{t} + \vK_{t}\delta \vx_{t}$. 
In  a nutshell, eq. \eqref{eq:state_DDP} becomes an equation which can be separated into constants and linear terms in $\delta \vx_{t}$. This is then used to generate the backward pass on $\bv_t (\bar{\blambda}_{t})$ and $\vS_{t}$. This process is repeated to obtain all $\delta \vu_{t}$ and update the trajectory. The general equations for backward pass are:
\begin{align}
    \bv_{t} = Q_{\vx}^{t} - (Q_{\vu}^{t})\tran (Q_{\vu\vu}^{t})^{-1}Q_{\vu\vx}^{t}, \label{eq:back_pass_DDP2}\\
    \vS_{t} = Q_{\vx\vx}^{t} - Q_{\vx\vu}^{t}(Q_{\vu\vu}^{t})^{-1}Q_{\vu\vx}^{t},
    \label{eq:back_pass_DDP}
\end{align}
where
\begin{align*}
    &Q_{\vx}^{t} = c_{\vx_{t}} +  \vf_{\vx_t}\tran \bv_{t+1} \\
    &Q_{\vu}^{t} = c_{\vu_{t}} +  \vf_{\vu_t} \tran \bv_{t+1}\\
    &Q_{\vx\vx}^{t} = c_{\vx_t\vx_t} + \vf_{\vx_t}\tran \vS_{t+1}\vf_{\vx_t} + \bv_{t+1}\otimes  \vf_{\vx_t \vx_t} \\
    &Q_{\vu\vu}^{t} = c_{\vu_t\vu_t} + \vf_{\vu_t}\tran \vS_{t+1}\vf_{\vu_t} + \bv_{t+1}\otimes \vf_{\vu_t \vu_t}\\
    &Q_{\vu\vx}^{t} = c_{\vu_t\vx_t} + \vf_{\vu_t}\tran \vS_{t+1}\vf_{\vx_t} + \bv_{t+1} \otimes \vf_{\vu_t \vx_t}
\end{align*}
with the terminal conditions:
\begin{align*}
    \vS_{T} = C_{T,\vx\vx} \quad \text{and} \quad \bv_{T} = C_{T,\vx}.
\end{align*}
The term $\vk_{t} = -(Q_{\vu\vu}^{t})^{-1}Q_{\vu}^{t}$ and $\vK_{t} = -(Q_{\vu\vu}^{t})^{-1}Q_{\vu\vx}^{t}$. DDP computes $\vk_{t}$ and $\vK_{t}$ for all time steps and then updates the control $\vu_{t} = \bar{\vu}_{t} + \alpha\vk_{t} + \vK_{t}\delta \vx_{t}$ where $\delta \vx_{t}$ is the perturbation given by the system after applying $\delta \vu_{t-1}$ which is $\delta \vx_{t} =  \vf(\bar{\vx}_{t-1} + \delta \vx_{t-1}, \bar{\vu}_{t-1} + \delta \vu_{t-1}) - \vf(\bar{\vx}_{t-1}, \bar{\vu}_{t-1})$, with $\delta \vx_{0} = 0$ as the initial state, $\vx_{0}$, is given.

\subsection{Stagewise Newton}

The method ``Stagewise Newton" solves the same set of necessary conditions as in eqs. \eqref{eq:lin_dyn}-\eqref{eq:con_DDP}, like DDP. However, there are two key differences. In the forward pass, it uses linearized dynamics to update the control. So, $\vu_{t} = \bar{\vu}_{t} + \alpha \vk_{t} + \vK_{t}\delta\hat{\vx}_{t}$, where $\delta \hat{\vx}_{t} = \vf_{\vx_{t-1}}\delta \hat{\vx}_{t-1} + \vf_{\vu_{t-1}}\delta \vu_{t-1}$. 

In addition to this, one replaces $\bar{\blambda}_{t} + \delta \blambda_{t} $ to a new variable $\vl_{t}$ in eqs. \eqref{eq:state_DDP} and \eqref{eq:con_DDP}. $\vl_{t}$ takes the form $\vS_{t}\delta \vx_{t} + \bv_{t}$ if the necessary conditions are solved in this manner. So, the necessary conditions become:
\begin{subequations}
    \begin{equation}
        \vf_{\vx_t}\delta \vx_t + \vf_{\vu_t}\delta \vu_t - \delta \vx_{t + 1} = \bo.
        \label{eq:lin_dyn_SN}
    \end{equation}
    \begin{equation*}
        c_{\vx_t} + c_{\vx_t \vx_t}\delta \vx_t + c_{\vx_t \vu_t}\delta \vu_t 
    \end{equation*}
    \begin{equation}
        +(\vf_{\vx_t \vx_t} \delta \vx_{t} + \vf_{\vx_t \vu_{t}}\delta \vu_{t})\tran \bar{\blambda}_{t+1} + \vf_{\vx_t}\tran \vl_{t+1}
        -\vl_t = 0
        \label{eq:state_SN}
    \end{equation}
    \begin{equation*}
        c_{\vu_t} + c_{\vu_t \vu_t}\delta \vu_t + c_{\vu_t \vx_t} \delta \vx_t
    \end{equation*}
    \begin{equation}
        +(\vf_{\vu_t \vu_t} \delta \vu_{t} + \vf_{\vu_t \vx_{t}}\delta \vx_{t})\tran \bar{\blambda}_{t+1} + \vf_{u_t}\tran \vl_{t+1} = 0
        \label{eq:con_SN}
    \end{equation}    
\end{subequations}
Additionally, rather than using $\bar{\blambda}_{t}$ from the previous iteration (which is done in a typical SQP solver), it enforces linearity of $\bar{\blambda}$ in terms of $\delta\vx$ and $\delta \vu$. The backward pass on $\bar{\blambda}$ comes from eq. \eqref{eq:state_SN}. Given that $\vl_{t}^{\text{iter}-1} = \bar{\blambda}_{t}^{\text{iter}}$ and $\vl_{t+1}^{\text{iter}-1} = \bar{\blambda}_{t+1}^{\text{iter}}$ , where ``iter" denotes the current iteration number. This along with the expansions
\begin{align*}
   &c_{\vx_t}^{\text{iter}} = c_{\vx_t}^{\text{iter}-1} + c_{\vx_t \vx_t}^{\text{iter}-1}\delta \vx_t + c_{\vx_t \vu_t}^{\text{iter}-1}\delta \vu_t  + \text{H.O.T.} \\
    &\vf_{\vx_t}^{\text{iter}} = \vf_{\vx_t}^{\text{iter}-1} + \vf_{\vx_t\vx_t}^{\text{iter}-1}\delta \vx_t + \vf_{\vx_t\vu_t}^{\text{iter}-1}\delta \vu_t + \text{H.O.T.}, 
\end{align*}
 and using eq. \eqref{eq:state_SN}, one can write
\begin{equation*}
    \bar{\blambda}_{t}^{\text{iter}} = \vl_{t}^{\text{iter}-1} =  \underbrace{c_{\vx_t} + c_{\vx_t \vx_t}\delta \vx_t + c_{\vx_t \vu_t}\delta \vu_t 
        +(\vf_{\vx_t\vx_t} \delta \vx_{t} + \vf_{\vx_t \vu_{t}}\delta \vu_{t})\tran \bar{\blambda}_{t+1} + \vf_{\vx_t}\tran \vl_{t+1}}_{\text{all derivatives at (iter - 1)}},     
\end{equation*}
which under linearity assumption gives:
\begin{equation}
    \bar{\blambda}_{t}^{\text{iter}} = c_{\vx_t}^{\text{iter}} +  (\vf_{\vx_t}^{\text{iter}})\tran \bar{\blambda}_{t+1}^{\text{iter}}.
    \label{eq:back_pass_SN_lam}
\end{equation}
However, this will have errors even under linearity assumption of derivatives as the term $(\vf_{\vx_t \vx_t} \delta \vx_{t} + \vf_{\vx_t \vu_{t}}\delta \vu_{t})\tran \bar{\blambda}_{t+1}$ should have been $(\vf_{\vx_t \vx_t} \delta \vx_{t} + \vf_{\vx_t \vu_{t}}\delta \vu_{t})\tran \vl_{t+1}$ for eq.~\eqref{eq:back_pass_SN_lam} to hold even in a linear sense. But, ``Stagewise Newton" makes this approximation.
So, the backward pass for ``Stagewise Newton" becomes 
\begin{align}
    \bv_{t} = Q_{\vx}^{t} - (Q_{\vu}^{t})\tran (Q_{\vu\vu}^{t})^{-1}Q_{\vu\vx}^{t}, \\
    \vS_{t} = Q_{\vx\vx}^{t} - Q_{\vx\vu}^{t}(Q_{\vu\vu}^{t})^{-1}Q_{\vu\vx}^{t},\\
    \bar{\blambda}_{t} = c_{\vx_t} + \vf_{\vx_t}\tran \bar{\blambda}_{t+1},
    \label{eq:SN_back_pass_lam}
\end{align}
where
\begin{align*}
    &Q_{\vx}^{t} = c_{\vx_{t}} + \vf_{\vx_t}\tran \bv_{t+1} \\
    &Q_{\vu}^{t} = c_{\vu_{t}} + \vf_{\vu_t}\tran \bv_{t+1}\\
    &Q_{\vx\vx}^{t} = c_{\vx_t\vx_t} + \vf_{\vx_t}\tran \vS_{t+1}\vf_{\vx_t} + \bar{\blambda}_{t+1}\otimes \vf_{\vx_t \vx_t}\\
    &Q_{\vu\vu}^{t} = c_{\vu_t\vu_t} + \vf_{\vu_t}\tran \vS_{t+1}\vf_{\vu_t} + \bar{\blambda}_{t+1} \otimes \vf_{\vu_t \vu_t}\\
    &Q_{\vu\vx}^{t} = c_{\vu_t\vx_t} + \vf_{\vu_t}\tran \vS_{t+1}\vf_{x_t} + \bar{\blambda}_{t+1} \otimes \vf_{\vu_t \vx_t}
\end{align*}
with the terminal conditions:
\begin{align*}
    \vS_{T} = C_{T,\vx\vx}, \quad \bv_{T} = C_{T,\vx}, \quad \text{and} \quad \bar{\blambda}_{T} = C_{T,\vx}. 
\end{align*}
We reiterate that the forward pass for Stagewise Newton updates $\delta \vu_{t} = \vk_{t} + \vK_{t} \delta \hat{\vx}_{t}$ where $\delta \hat{\vx}_{t} = \vf_{\vx_{t-1}}\delta \hat{\vx}_{t-1} + \vf_{\vu_{t-1}}\delta \vu_{t-1}$, is the linearized dynamics, unlike DDP, which takes the full $\delta \vx_{t} = \vf(\bar{\vx}_{t-1} + \delta \vx_{t-1},\bar{\vu}_{t-1}+\delta \vu_{t-1}) - \vf(\bar{\vx}_{t-1},\bar{\vu}_{t-1})$.
\subsection{iLQR}
We describe the iLQR algorithm next. iLQR uses the modified QP as mentioned in eq. \eqref{eq:iLQR-LQR}. For the optimal control problem, this means that the optimization problem in eq. \eqref{eq:per_optim} becomes
\begin{dmath}
    \min_{\delta \vx_t, \delta \vu_{t}, \delta \blambda_t} \sum_{t=0}^{T-1}\Big(c_{\vx_t}\delta \vx_{t} + c_{\vu_t}\delta \vu_{t} + \frac{1}{2}\delta \vx_{t}' c_{\vx_t\vx_t} \delta \vx_t + \delta \vx_t\tran c_{\vx_t\vu_t} \delta \vu_t + \frac{1}{2} \delta \vu_{t}\tran c_{\vu_t \vu_t} \delta \vu_t  + (\bar{\blambda}_{t+1} + \delta \blambda_{t+1})\tran (\vf(\bar{\vx}_t, \bar{\vu}_t) + \vf_{\vx_t}\delta \vx_t + \vf_{\vu_t} \delta \vu_t - 
    \bar{\vx}_{t+1}  -\delta \vx_{t+1})\Big ) + C_{T, \vx}\delta \vx_T + \frac{1}{2} \delta \vx_T\tran C_{T, \vx\vx} \delta \vx_T.
    \label{eq:per_optim_iLQR}
\end{dmath}
It has the same forward and backward pass as DDP except we drop the tensor products of second order terms in dynamics. This gives:

\begin{align*}
    &Q_{\vx}^{t} = c_{\vx_{t}} + \vf_{\vx_t}\tran\bv_{t+1} \\
    &Q_{\vu}^{t} = c_{\vu_{t}} + \vf_{\vu_t}\tran \bv_{t+1}\\
    &Q_{\vx\vx}^{t} = c_{\vx_t\vx_t} + \vf_{\vx_t}\tran \vS_{t+1}\vf_{\vx_t} \\
    &Q_{\vu\vu}^{t} = c_{\vu_t\vu_t} + \vf_{\vu_t}\tran \vS_{t+1}\vf_{\vu_t} \\
    &Q_{\vu\vx}^{t} = c_{\vx_t\vu_t} + \vf_{\vu_t}\tran \vS_{t+1}\vf_{\vx_t} 
\end{align*}
Everything else follows the same as DDP.

\end{document}